\documentclass[journal, draftcls]{IEEEtran}

\onecolumn

\newcommand\nval{10,000}
\newcommand\nMonte{10,000}

\newcommand{\qed}{\hfill \ensuremath{\Box}}

\RequirePackage{amsthm,amsmath,amsfonts,amssymb}
\RequirePackage[colorlinks,citecolor=blue,urlcolor=blue]{hyperref}
\RequirePackage{graphicx}
\usepackage{nccmath}
\usepackage{subfiles}
\usepackage{tikz}
\usepackage[most]{tcolorbox}
\usepackage{pgfplots}
\pgfplotsset{compat=1.18}
\usepackage{color}
\usepackage{subfiles}

\newtheorem{theorem}{Theorem}
\newtheorem{prop}[theorem]{Proposition}
\newtheorem{cor}[theorem]{Corollary}
\newtheorem{lem}[theorem]{Lemma}

\newcommand{\col}{\mathsf{col}}

\newcommand{\Pois}{\mathrm{Pois}}

\newcommand{\Po}{\mathrm{P}}

\newcommand{\Prp}[1]{\Pr\left[#1 \right]}
\newcommand{\reals}{\mathbb R}

\newcommand{\Ncal}{\mathcal{N}}

\newcommand{\one}{\mathbf{1}}

\newcommand{\Var}[1]{\mathrm{Var}\left[ #1\right]}
\newcommand{\ex}[1]{\ensuremath{\mathbb{E}\left[ #1\right]}}
\newcommand{\exsub}[2]{\ensuremath{\mathbb{E}_{#1}\left[ #2\right]}}

\newcommand{\notered}[1]{{\sf\textcolor{red}{[#1]}}}

\newcommand{\Cov}{\mathrm{Cov}}

\newcommand{\diag}{\operatorname{diag}}

\begin{document}

\author{Alon~Kipnis,~\IEEEmembership{Member,~IEEE,}
\thanks{This work is partially funded by the US-Israel Binational Science Foundation under grant number 2022124.}
\thanks{Alon Kipnis is with the School of Computer Science at Reichman University, Herzliya, Israel.}
\thanks{This paper was presented in part at the 2023 59th Annual Allerton Conference on Communication, Control, and Computing 
\cite{KipnisAllerton2023minimax}. }
}

\title{The Minimax Risk in Testing Uniformity over Large Alphabets under Missing-Ball Alternatives}


\maketitle

\begin{abstract}
We study the problem of testing the goodness of fit of categorical count data to a Poisson distribution uniform over the categories, against a class of alternatives defined by excluding an $\ell_p$ ball, $p \leq 2$, of radius $\epsilon$ around the uniform rate sequence. We characterize the minimax risk for this problem as the expected number of samples $n$ and the number of categories $N$ go to infinity. Our result enables constant-factor comparisons among the many estimators previously proposed for this problem, rather than comparisons only at the level of convergence rates or scaling orders of sample complexity. The minimax test relies exclusively on collisions in the small sample limit, but behaves like the chi-squared test otherwise. Empirical studies across a range of parameters show that the asymptotic risk estimate is accurate in finite samples, and that the minimax test outperforms both the chi-squared test and a test based on collisions under the least favorable alternative. Our analysis involves a reduction to a structured subset of alternatives, establishing uniform asymptotic normality for a family of linear test statistics, and solving an optimization problem over $N$-dimensional sequences akin to classical results from signal detection in Gaussian white noise. Finally, we discuss the connection to the fixed-sample-size multinomial model, arguing that the Poisson minimax risk derived here also characterizes the minimax risk of the multinomial problem.

\end{abstract}

\section{Introduction}
\label{sec:intro}

\subsection{Background and Motivation}

We observe data consisting of occurrence counts 
$O_1,\ldots,O_N$, where each $O_i$ corresponds to one of $N$ distinct categories. We model these counts as realizations of independent Poisson random variables and are interested in testing whether all categories share the same Poisson rate, or whether the rates vary across categories according to some alternative structure. Specifically, suppose
\begin{align}
\label{eq:hyp_Q}
H(Q)\,:\,O_i \sim \Pois(nQ_i), 
\end{align}
independently for $i=1,...,N$, where $Q = (Q_1,\ldots,Q_N) \in \reals_+^N$ is an arbitrary Poisson rate vector and $n$ is a scaling parameter that controls the total expected number of occurrences. We are interested in testing the null hypothesis 
\begin{align*}
H_0: Q_i=1/N \quad \text{for all}
\quad i=1,\ldots,N,
\end{align*}
against alternatives in which the rates $Q_i$ deviate from uniformity. Our results are asymptotic and valid across all scaling relationships between $N$ and $n$. Still, we focus on the high-dimensional regime ($N \gg n$), which is of particular interest, as the case where $N$ is comparable to or smaller than $n$ is well understood.

This situation might arise when we are interested in testing the spatial uniformity of a light source using a large array of photon counters operating over a prescribed time interval, or when verifying whether the symbols in a random string are distributed uniformly over a large alphabet when the string's length has a Poisson distribution with mean $n$. In biology, the counts might represent the number of observed individuals of a particular plant or animal species in $N$ different quadrats within a larger habitat. Testing whether these counts follow the same Poisson law would help assess if the species is uniformly distributed across the habitat or if there are preferred sub-regions due to variations in environmental conditions, leading to inhomogeneous Poisson rates \cite{illian2008statistical}.


A fundamental question regarding the Poisson model is its relationship to the multinomial setting, which is often the primary interest in applications. Conditioning the vector of independent Poisson counts in \eqref{eq:hyp_Q} on the event $\sum_{i=1}^N O_i = n$ yields the multinomial distribution with cell probabilities proportional to $Q$, provided $\sum_{i=1}^N Q_i = 1$ \cite[Ch. 6.3]{barbour1992poisson}. Consequently, the minimax risk for the Poisson model bounds from above the minimax risk in the multinomial model.
In fact, the technique developed in this paper extends beyond this upper-bound relationship: the minimax asymptotic derived here also characterizes the corresponding quantity under multinomial sampling. The extension requires additional technical steps to address the dependence introduced by conditioning, which can be handled using suitable conditional central limit theorems for count data \cite{janson2001moment,gnedin2007notes,klein2015conditional}. A detailed treatment of the multinomial case, establishing this equivalence, is presented in a companion note \cite{kipnis2026CLT}. These connections align our results with several classical Poissonization/de-Poissonization equivalence theorems \cite{renyi1962three,holst1972asymptotic,morris1975central,quine1984normal,Nussbaum1996,janson2001moment}. However, the high-dimensional regime and the nature of the alternatives considered here introduce new challenges not covered by existing results; see Section~\ref{sec:discussion} for a detailed discussion.


Testing for uniformity includes, as a special case, testing the goodness-of-fit to any continuous distribution with cumulative distribution function (CDF) $F_0$. This is done by reducing the sample $y_1,\ldots,y_m$ to bin counts
\begin{align}
    O_i = \#\left\{ j \in \{1,\ldots,m\}\,:\,\frac{i-1}{N} \leq F_0(y_j) < \frac{i}{N} \right\}, \quad i=1,\ldots,N;
    \label{eq:uniform}
\end{align}
as discussed in \cite{ermakov1998asymptotic} and \cite[Ch. 1.4]{ingster2003nonparametric}. This focus on continuous distributions is not restrictive in high dimensions (large $N$), where practical constraints often necessitate assuming a smooth, low-dimensional model due to the difficulty of estimating complex distributions from limited data.

\subsection{Problem Formulation}
Denote by $U = (1/N,...,1/N) \in \reals_+^N$ the uniform rate vector. Given a test $\psi : \mathbb N^N \to \{0,1\}$ and $Q \neq U$, the risk of $\psi$ is 
\begin{align*}
    R(\psi; Q) & := \Prp{\psi(O_1,\ldots,O_N) \mid H_0} \\
    & \qquad + \Prp{1-\psi(O_1,\ldots,O_N) \mid H(Q)}.  
\end{align*}
Given a non-empty set of alternative sequences $V$, we are interested in the minimax risk over $V$:
\begin{align}
    \label{eq:minimax_risk}
    R^*(V) := \inf_{\psi} \sup_{Q \in V} R(\psi; Q). 
\end{align}
We also use the convention $R^*(\emptyset)=0$, which is natural by the previous definition.  

The minimax analysis is commonly understood as a two-person game of the statistician versus Nature \cite[Ch. 5]{berger2013statistical},\cite[Ch. 6]{gine2021mathematical}: The statistician plays an estimator $\psi$ to decide whether the data generating frequency sequence $Q\in \reals^N_+$ is in the null or the alternative. Nature plays a choice of $Q$, either from the null or the alternative. Nature tries to maximize the risk while the statistician tries to minimize it. 

In this paper, we consider a set of alternative rate sequences $V_\epsilon$  obtained by removing an $\ell_p$ ball of radius $\epsilon$ around the uniform rate sequence $U$, for $p \leq 2$. Namely,
\begin{align}
    \label{eq:alternative_def}
V_{\epsilon} := \left\{ Q \in \reals_+^N\,:\, 
\left\| Q -U \right\|_p \geq \epsilon
\right\},\qquad p \in (0,2],
\end{align}
where $\left\|a \right\|_p = \left(\sum_{i=1}^N \left|a_i\right|^p \right)^{1/p}$ is the $\ell_p$ norm in $\reals^N$. 

\begin{figure}
    \centering
    \begin{tikzpicture}
    \def\r{2}
    \def\cc{\r*0.86}
    \node[coordinate] (center) at (0,0) {};
    \node[coordinate] (lf) at (\r/1.55,\r/1.55) {};
 
    \draw[fill=red, fill opacity=0.25, draw=white] (-\cc,-\cc) rectangle (1.5*\cc,1.5*\cc);


    \draw[fill=white] (-\r,0)  to[bend left] (0,\r)  to[bend left] (\r,0)  to[bend left] (0,-\r) to[bend left] (-\r,0);

    \draw[<->,color=red, dotted] (0,0) -- node[above] {$\mu^*$} (lf);

    \draw[fill=red, color=red] (\r/1.55,\r/1.55) circle (1pt);
    
    \draw[fill=red, color=red] (\r/1.55,-\r/1.55) circle (1pt);
    \draw[fill=red, color=red] (-\r/1.55,-\r/1.55) circle (1pt);
    \draw[fill=red, color=red] (-\r/1.55,\r/1.55) circle (1pt);
    \draw[<->] (0,0) -- node[above] {$\epsilon$} (\r,0);
    \draw[fill=blue] (center) circle (1pt) node[below, color=blue, xshift=0cm] {\footnotesize $\left(\frac{1}{N},\frac{1}{N}\right)$};
    
    \node[below left] at (1.2*\cc,1.2*\cc) {\small $V_\epsilon$};
    
    \draw (-\r,0) 
            to[bend left] (0,\r) 
            to[bend left] (\r,0) 
            to[bend left] (0,-\r) 
            to[bend left] (-\r,0);
    \node (northwest) at (-0.85*\cc,0.88*\cc) {};
    \node (southwest) at (-0.84*\cc,-0.88*\cc) {};
    \node (southeast) at (0.85*\cc,-0.9*\cc) {};
    \node (northeast) at (0.85*\cc,0.89*\cc) {};
    
    \draw[->, color=black, dotted] (-\cc,-
    \cc) -- (-\cc,1.6*\cc);
    \draw[->, color=black, dotted] (-\cc,-\cc) -- (1.7*\cc,-\cc);
    \draw[fill=black, color=black] (-\cc,-\cc) circle (1pt) node[below] {\footnotesize $(0,0)$};
    
    \end{tikzpicture}
    \caption{Conceptual sketch of the sets of alternatives $V_{\epsilon}$ (shaded red) in $N=2$ dimensions and some $p \in (1,2]$. The least favorable rate sequences in $V_\epsilon$ are typical realizations of a prior supported by the points at the boundary of the $\ell_p$ ball around the uniform rate sequence $U = (1/N,...,1/N)$ closest to the center (indicated by 4 red dots). $\mu^*=\epsilon N^{-1/p}$ is the perturbation defining the least favorable prior.}
    \label{fig:alternative}
\end{figure}

Throughout this paper, we use the standard $o$ and $O$ notations to denote asymptotic relations between sequences of real numbers. For example, $f(N) = o(g(N))$ means that $\lim_{N \to \infty} f(N)/g(N) = 0$. 

\subsection{Previous Work}
The minimax risk in the case $N = o(n)$ is well-known and follows from the minimaxity of the chi-squared test (c.f. \cite[Ch. 1]{ingster2003nonparametric}). The focus of this paper is the case $n = o(N)$. This setting is related to non-parametric hypothesis testing on densities \cite{ingster1987minimax,ermakov1990asymptotically,ingster1993asymptotically,lepski1999minimax,balakrishnan2019hypothesis} and to testing for the uniformity under the multinomial model mentioned earlier \cite{diaconis1989methods,batu2001testing,paninski2008coincidence,balakrishnan2018hypothesis}. 
In these contexts, \cite{ingster1993asymptotically,ingster2003nonparametric} characterized the minimax risk when each sequence in the alternative is a binned version of a smooth density function as in \eqref{eq:uniform} and showed that the minimax test is based on the chi-squared statistic. Similar results under assumptions other than smoothness and some different alternatives can be found in \cite{ermakov1998asymptotic}, including asymptotic minimaxity of the chi-squared test under $\ell_2$ alternatives. In recent years, works originating from the field of property testing in computer science \cite{ron2008property,goldreich_2017} focused on testing uniformity against discrete distribution alternatives that do not necessarily arise as binned versions of smooth densities
\cite{batu2001testing,paninski2008coincidence,diakonikolas2016new,acharya2018improved,waggoner2015lp}. Instead, they may be unrestricted or obey other properties \cite{rubinfeld2005testing,diakonikolas2014testing,gs009}, and typically focus on the case $p=1$. These works characterized estimators' optimal rate of convergence, e.g., the number of samples guaranteeing vanishing minimax risk in the other problem parameters. Nevertheless, these previous works neither provide the asymptotic minimax risk nor identify the minimax test in either the Poisson or multinomial setting, which have remained open problems. The present work delivers the minimax risk in both settings.

\subsection{Contributions}
Consider an asymptotic setting where $N$ and $n$ go to infinity. If $\lim_{n,N \to \infty} R^*(V_\epsilon) \in (0,1)$, then 
\begin{align}
\label{eq:result_intro}
 \lim_{n,N \to \infty} \frac{R^*(V_\epsilon)}{2 \Phi \left(- u_{\epsilon,n,N,p}/2 \right)} =  1.
\end{align}
where
 \begin{align}
(u_{\epsilon,n,N,p})^2:= \frac{1}{2}\epsilon^4 n^2 N^{3-4/p}. 
\label{eq:u_def} 
\end{align}
Furthermore, sufficient conditions for $\lim_{N,n \to \infty} R^*(V_\epsilon) \in (0,1)$ are 
$u_{\epsilon,n,N,p}\to c$ and $N=o(n^2)$ or $\epsilon N^{1-1/p} = o(1)$. Under these conditions, a test statistic linear in the histogram ordinates 
\begin{align}
    \label{eq:Xm_def}
X_m := \sum_{i=1}^N \one\{O_i=m\},\qquad m=0,1,2,\ldots.
\end{align}
is asymptotically minimax, and the least favorable prior in a Bayesian counterpart of the minimax problem is an $N$ product of a symmetric two-point prior with support $\left\{1/N - \epsilon N^{-1/p}, 1/N + \epsilon N^{-1/p}\right\}$; see Figure~\ref{fig:alternative} for a conceptual sketch. 

Additionally, we derive the asymptotic risk of the chi-squared test under the least favorable prior. Under the same conditions as above, this risk converges to
\begin{align}
\label{eq:chi2_intro}
2 \Phi\left(-\sqrt{\frac{n}{N+n}}\frac{u_{\epsilon,n,N,p}}{2}\right).
\end{align}
This expression shows that, unless $n/N \to \infty$, the chi-squared test fails to achieve the minimax risk under Poisson sampling. We emphasize that this phenomenon is not merely a consequence of an incorrect null calibration: it is well known that the classical chi-squared statistic does not follow a chi-squared distribution under the null when $n/N \to 0$ \cite{haberman1988warning}. In contrast, the test analyzed here employs an optimally chosen threshold for the chi-squared statistic so as to minimize the total risk.

Numerical analyses in the case $p=1$ show that the approximation \eqref{eq:result_intro} provides to $R^*(V_\epsilon)$ is accurate in finite $n$, $N$, and small values of $\epsilon$; see for example Figure~\ref{fig:risk} below. These analyses also demonstrate the dominance of the minimax test over the chi-squared test and a test based on collisions \cite{diakonikolas2019collision}. 

\subsection{Significance of the contributions}
Previous results in the literature showed that the condition $u_{\epsilon,n,N,p} \to \infty$ implies complete separation ($R^*(V_\epsilon) \to 0$) and $u_{\epsilon,n,N,p} \to 0$ implies the impossibility of separation ($R^*(V_\epsilon) \to 1$) \cite{paninski2008coincidence,ValiantValiant,chhor2022sharp}. These results established the rate optimality of the problem. Our main results characterize the asymptotic risk in the entire regime $R^*(V_\epsilon) \to c \in (0,1)$. 
Our characterization enables a principled comparison among various estimators that achieve the minimax sample complexity, including those studied in \cite{batu2001testing,paninski2008coincidence,chan2014optimal,diakonikolas2016new,balakrishnan2019hypothesis,chhor2022sharp,gupta2022sharp}. The minimax risk captures performance at the level of leading constants, akin to the role of Pinsker's constant in nonparametric function estimation and Fisher information in parametric inference; see the discussion in \cite{nussbaum1999minimax}.

The minimax properties of the chi-squared test serve as a particularly illustrative case. It follows from previous works in the multinomial that the chi-squared test attains vanishing minimax risk as $u_{\epsilon,n,N,1} \to \infty$, thus it is order-optimal in the case $p=1$ \cite{acharya2015optimal,chhor2022sharp}. Furthermore, the chi-squared test is asymptotically minimax either when the alternatives are smooth \cite{ingster1993asymptotically} or when $p=2$, $N = o(n^2)$, and $\epsilon \sim 1/\sqrt{N}$ \cite{ermakov1998asymptotic}. The findings of this paper, in particular equations \eqref{eq:result_intro} and \eqref{eq:chi2_intro}, precisely characterize the extent to which the chi-squared test deviates from the minimax optimal test under Poisson sampling. 
%
%

\subsection{Proof Technique}
We derive both an upper bound (Theorem~\ref{thm:upper_bound}) and a lower bound (Theorem~\ref{thm:lower_bound}) for $R^*(V_{\epsilon})$, and show that these bounds are asymptotically tight whenever $R^*(V_\epsilon)$ does not vanish asymptotically (Theorem~\ref{thm:matching}). Both bounds are established via a related Bayesian framework in which Poisson rate sequences are drawn from a class of priors associated with the most informative subset of $V_\epsilon$. This technique is standard in minimax analysis when direct characterization of the least favorable prior is intractable \cite{donoho1994minimax,ingster1993asymptotically,suslina1999extremum,lepski1999minimax,donoho2006asymptotic}. However, our setting does not reduce to any of these classical frameworks, as the counts, or their histogram ordinates, are generally not asymptotically normal \cite{kimber1983note}. The main technical difficulty lies in the derivation of the upper bound. This involves reducing the analysis to a carefully constructed subset of alternative rate vectors, defined via simple separation conditions; establishing an asymptotic normality result for a class of tests linear in the counts histogram; identifying suitable priors over the reduced alternative set; and ultimately solving an optimization problem over $N$-dimensional sequences. The structure of this optimization problem resembles those arising in hypothesis testing within the Gaussian white noise model \cite{ingster1993asymptotically,suslina1999extremum}.

An additional complication arises from alternatives that do not give rise to asymptotically normal test statistics. These are handled in Section~\ref{sec:prep} using standard tests whose behavior can be analyzed via their first two moments. 

\subsection{Paper outline}
Preliminary results are provided in Section~\ref{sec:prep}. We present and discuss the main results in Section~\ref{sec:main_results}. In Section~\ref{sec:numerical_evalutations}, we report on numerical simulations. Additional discussion and remarks are in Section~\ref{sec:discussion}. All the proofs are in Section~\ref{sec:proofs}, with some more technical details deferred to the appendix. 

\section{ Preliminaries \label{sec:prep}}
In this section, we develop several intermediate results and technical tools that are used in the proofs of the main results. 

\subsection{Elementary Separation Conditions} 

We first examine several test statistics used to identify subsets of $V_\epsilon$ where the minimax risk vanishes. These subsets correspond to alternatives that are easily distinguishable from the null and are subsequently excluded from the main analysis, allowing us to focus on the more challenging regions of $V_\epsilon$.

For a test statistic $T$ and an alternative rate sequence $Q$, define 
\[
A_{n,N}(T;Q) := \frac{\ex{T \mid H(Q)} - \ex{T \mid H(U)} }{\sqrt{ \Var{T \mid H(Q)} + \Var{T \mid H(U)}}}. 
\]
When $A_{n,N}(T;Q) \to \infty$, one can construct a test based on $T$ with vanishing risk using a standard argument via Chebyshev's inequality. For example, reject the null when 
\[
T \geq \ex{T \mid H(U)} + \frac{A_{n,N}(T;Q)}{2} \sqrt{\Var{T \mid H(U)}}.  
\]
In this case, the null $H(U)$ and the alternative $H(Q)$ are said to be completely separated. 

\subsubsection*{Sum Test}
Consider the statistic
\[
T_{\mathrm{sum}} := \sum_{i=1}^N O_i. 
\]
Direct calculation shows 
\begin{align*}
A_{n,N}(T_{\mathrm{sum}};Q) & = \frac{n \sum_{i=1}^N Q_i - n}{\sqrt{n + n \sum_{i=1}^N Q_i }}. 
\end{align*}
It follows that a test based on $T_{\mathrm{sum}}$ can completely separate an alternative $Q \in V_\epsilon$ unless 
\begin{equation}
\label{eq:trivial_sep_cond}
    \sqrt{n}\left(\sum_{i=1}^N Q_i - 1\right) = O(1).
\end{equation}

\subsubsection*{Chi-squared Test}
The chi-squared test is based on the statistic
\begin{align}
    \label{eq:chisquared_def}
T_{\chi^2} := \sum_{i=1}^N \frac{(O_i - n/N)^2}{n/N}. 
\end{align}
This test is known to yield an asymptotically minimax test when $N = o(n)$ or when the set of alternative sequences is restricted to a binned smooth probability density \cite{ingster1987minimax,ermakov1998asymptotic,ingster2003nonparametric}. 

If $Q$ satisfies \eqref{eq:trivial_sep_cond}, then the variance of $T_{\chi^2}$ under $H(Q)$ is at most $3 N(1 + O(1))$. In this case,
\begin{align*}
A_{n,N}(T_{\chi^2};Q) & \geq \frac{-1+N \sum_{i=1}^N Q_i + \frac{N}{n} \sum_{i=1}^N \left(n Q_i- n/N \right)^2 }{2\sqrt{3N(1+O(1))}} \\
& \geq C n \sqrt{N} \left\|Q - U\right\|_2^2 + o(1). 
\end{align*}
for some constant $C>0$. Therefore, the chi-squared test completely separates alternatives outside the ball $\left\|Q - U\right\|_2 \leq r_{n,N}$, provided
\begin{align}
    \label{eq:L2_sep_condition}
    r_{n,N}^2 n \sqrt{N} \to \infty. 
\end{align}
Furthermore, using the inequality $\|a\|_p \leq N^{1/p-1/2} \|a\|_2$ valid for all $a \in \reals^N$ and $p\leq 2$, we get complete separation whenever 
\begin{align}
    n N^{3/2-2/p} \left( \left\|Q - U \right\|_p \right)^2  \to \infty. 
    \label{eq:trivial_separation_cond_chisq_general}
\end{align}
For $Q \in V_\epsilon$, condition \eqref{eq:trivial_separation_cond_chisq_general} is equivalent to $u_{\epsilon,n,N,p} \to \infty$. This condition was recognized as sufficient for complete separation in \cite{chhor2022sharp} and in \cite{paninski2008coincidence,ValiantValiant,gupta2022sharp} under the related multinomial setups. These works also showed a converse statement: $u_{\epsilon,n,N,p} \to 0$ leads to inseparability, i.e. minimax risk converging to one.

\subsubsection*{Max Test}
Consider the statistic 
\[
T_{\max} :=  \max_{i=1,\ldots,N} O_i. 
\]
Under the null $Q = U$, the distribution of $T_{\max}$ is degenerate and concentrates on two consecutive integers approaching $\log(N)/\log(\log(N))$ to first order in $N$ \cite{kolchin1969limiting,kimber1983note}. Therefore, under any $Q\neq U$, if $nQ_i \to \infty$ faster than $\log(N)/\log(\log(N))$ for some $i$, a test that rejects the null if $T_{\max}$ exceeds $2 + \log(N)/\log(\log(N))$ detects with probability of error approaching zero and thus has vanishing risk. Consequently, a test based on $T_{\max}$ completely separates alternatives outside the hypercube 
\begin{align}
    B_\xi^\infty(U) := \{Q \in \reals_+^N\,:\, \max_{i=1,\ldots,N} |Q_i-U_i| \leq \xi_{n,N}\},
    \label{eq:hypercube_def}
\end{align}
provided
\begin{align}
 \frac{n\log(\log(N))}{\log(N)} \xi_{n,N} \to \infty.
\label{eq:max_test_cond}
\end{align}

The following corollary summarizes the separation of the chi-squared and max tests in our setting.  
\begin{cor}
\label{cor:trivial_separation}
Consider a sequence of multivariate Poisson models \eqref{eq:hyp_Q} indexed by $n$ and $N$, where $N$ and $n$ go to infinity. 
    Let $\xi = \xi_{n,N}$ satisfy \eqref{eq:max_test_cond}. Then
    \[
 R^*\!\left( V_{\epsilon} \cap B_\xi^\infty(U) \right) = R^*(V_{\epsilon}) + o(1). 
    \]
    Additionally, let $r=r_{n,N}$ satisfy \eqref{eq:L2_sep_condition}. Then
       \[
 R^*\!\left(V_{\epsilon} \cap B_\xi^\infty(U)  \cap B_{r}^2(U) \right) = R^*(V_{\epsilon}) + o(1). 
    \]
\end{cor}


\subsection{Minimax Bayesian setup}
Assume that the sequence $Q$ is sampled from some prior $\pi$ over $\reals_+^N$. The Bayes risk of a test $\psi$ is defined as
\begin{align*}
    \rho(\psi;\pi) =  \exsub{Q \sim \pi}{R(\psi; Q)},
\end{align*}
where we used the notation 
\begin{align*}
\exsub{Q\sim \pi}{F(Q)} & = 
 \int_{\reals^N} F(q_1,...,q_N) \prod_{i=1}^N \pi_i(dq_i) 
\end{align*}
for a measurable function $F : \reals^N \to \reals$, assuming all per-coordinate integrals exist. 

Consider the restricted set of Poisson rate sequences 
\begin{align}
    \label{eq:V_eps_xi_r_def}
V_{\epsilon,\xi,r} := V_\epsilon \cap B^\infty_\xi(U) \cap B^2_r(U).
\end{align}
By the inequality $\|a\|_p \leq N^{1/p-1/2} \|a\|_2$ mentioned earlier, $V_{\epsilon,\xi,r}$ is non-empty provided
\begin{align}
         \epsilon N^{1-1/p} \leq \xi \quad \text{and} \quad \epsilon N^{1/2-1/p} \leq r.
\label{eq:non_empty_constraint}
\end{align}
We consider a companion set to $V_{\epsilon,\xi,r}$ of product priors $\pi$ over $\reals^N$:
\begin{align}
\label{eq:Pi_eps_def}
    \Pi_{\epsilon,\xi,r}  := \Biggl\{ \pi = \prod_{i=1}^N \pi_i \,:& \, \sum_{i=1}^N \exsub{Q_i \sim \pi_i}{\left|Q_i-1/N \right|^{p}} \geq \epsilon^p,  \\
    & \left. \sum_{i=1}^N \exsub{Q_i \sim \pi_i}{\left|Q_i-1/N \right|^{2}} \leq r^2 \right.,  \nonumber \\
    & \pi([1/N-\xi, 1/N+\xi])=1 \Biggr\}. \nonumber
\end{align}
The minimax Bayes risk over $\Pi_{\epsilon,\xi,r}$ is defined as
\[
\rho^*(\Pi_{\epsilon,\xi,r}) := \inf_{\psi} \sup_{\pi \in \Pi_{\epsilon,\xi,r}} \rho(\psi;\pi).
\]
Whenever it exists, a prior $\pi^* \in \Pi_{\epsilon,\xi,r}$ attaining the supremum above is said to be \emph{least favorable}. 
\begin{lem}
\label{lem:R_geq_rho}
     We have
    \begin{align}
\label{eq:equality_R_rho}
    R^*(V_{\epsilon,\xi,r}) = \rho^*(\Pi_{\epsilon,\xi,r}) + o(1). 
    \end{align}
\end{lem}
The proof of Lemma~\ref{lem:R_geq_rho} is given in Section~\ref{sec:R_geq_rho:proof}. In \eqref{eq:equality_R_rho}, $o(1)$ represents a term that goes to zero as $N$ goes to infinity independently of the other parameters. The situation of most interest to us is when the risk on either side of \eqref{eq:equality_R_rho}
converges to a positive constant, in which case 
\eqref{eq:equality_R_rho} implies $R^*(V_{\epsilon,\xi,r})/\rho^*(\Pi_{\epsilon,\xi,r}) \to 1$. 

\subsubsection*{Example: Three-point prior}
As an example of an interesting set of priors, consider a prior in which each coordinate $\pi_i$ is a three-point (or two-point if $\eta=1$) univariate prior symmetric around $U_i=1/N$: 
\begin{align}
\pi_i(\eta,\mu) = (1-\eta)\delta_{U_i} + \frac{\eta}{2} \delta_{U_i+\mu} + \frac{\eta}{2} \delta_{U_i-\mu}, 
\label{eq:three_point_prior}
\end{align}
for $i=1,\ldots,N$. One of our key results says that the asymptotically unique least favorable prior $\pi^*$ within $\Pi_{\epsilon,\xi,r}$ is of the form \eqref{eq:three_point_prior}. 

\subsection{Properties of the counts' histogram}
\label{sec:histogram}
Recall that $X_m$ of \eqref{eq:Xm_def} denotes the number of categories with exactly $m$ items. For example, $X_0$ is the number of categories not represented in the data, $X_1$ is the number of singletons, and $X_2$ is the number of exclusively colliding pairs. We call the set
$\{X_0,X_1,\ldots\}$ the data's histogram 
(aka the data's \emph{pattern}  \cite{acharya2009maximum} or \emph{fingerprint} \cite{valiant2017estimating}). 

For $\lambda >0$, let $\Po_{\lambda}(m)$ be the Poisson probability mass function:
\[
\Po_{\lambda}(m) := e^{-\lambda}\frac{\lambda^m}{m!},\quad m=0,1,...
\]
Set
\[
\lambda_0:= \lambda_{n,N} := nU_i = n/N.
\]
Under the null, the histogram ordinate $X_m$ is a binomial random variable with mean 
\begin{align*}
    \mu^0_m &:= \ex{X_m} = \sum_{i=1}^N \Po_{\lambda_0}(m) = N \cdot \Po_{\lambda_0}(m).
\end{align*}
The covariance between the $\{X_m\}$s is
\begin{align}
    \Cov[X_m, X_k] = N\begin{cases}
        \Po_{\lambda_0}(m)(1 - \Po_{\lambda_0}(m)) & m=k \\
        - \Po_{\lambda_0}(m)\Po_{\lambda_0}(k) & m \neq k. 
    \end{cases} 
    \label{eq:cov}
\end{align}
Henceforth, we write the covariance function using the infinite matrix $\Sigma$ such that $\Sigma_{m,k} = \Cov[X_m, X_k]$, where we agree that the first row and column of $\Sigma$ have index $0$. It is convenient to write 
\begin{align}
    \label{eq:Sigma_matrix}
\Sigma = \diag(\mu^0) - \mu^0 \mu^{0\top}/N,
\end{align}
where for two sequences $u=\{u_m\}_{m=0}^\infty$ and $v=\{v_m\}_{m=0}^\infty$ the notation $u v^\top$ denotes the infinite matrix whose $(m,k)$ entry is $u_m v_k$, in accordance with standard matrix notation. We note that $\Sigma$ is singular because $\Sigma \one =0$, where $\one = (1,1,\ldots)$.

\subsection{Linear tests of the histogram}
\label{sec:analysis:linear_tests}

Consider tests that reject the null for large values of 
\begin{align}
    \label{eq:linear_test_stat}
T_w = \left\langle w, X \right\rangle =\sum_{m=0}^\infty X_m w_m 
\end{align}
for some weights sequence $ \{w_m\}_{m=0}^\infty$. Note that
\begin{align}
    \label{eq:T_as_sum_of_w}
T_w = \sum_{i=1}^N A_i,\qquad A_i := w_{O_i}, 
\end{align}
and that $A_i$, $i=1,2,...$ are independently and identically distributed under the null. Under additional assumptions on the sequence $\{w_m\}$  provided in Proposition~\ref{prop:asymp_normality} below, $T_w$ is asymptotically normal with mean and variance
\begin{align*}
    \ex{T_w} & = \sum_{m=0}^\infty w_m \mu^0_m =  \langle w, \mu^0 \rangle, \\
    \Var{T_w} & = \langle w, \Sigma w \rangle. 
\end{align*}
In this situation, a test asymptotically of size $\alpha$ against $H_0$ is obtained by rejecting when 
\begin{align}
    \label{eq:test_alpha_general}
\frac{T_w - \langle w, \mu^0 \rangle}{\sqrt{\langle w, \Sigma w \rangle}} > z_{1-\alpha},\qquad \Phi(z_{1-\alpha})=1-\alpha.  
\end{align}
In the sequel, we use tests of the form \eqref{eq:test_alpha_general} with $\alpha$ optimized to minimize the risk (Proposition~\ref{prop:asymp_normality}-(ii) below). 

We now analyze the mean shift in the statistic $T_w$ of the form \eqref{eq:linear_test_stat} under some prior $\pi$ for $Q$. For $x \in \reals$, $\lambda>0$, and $m=0,1,\ldots$, define
\begin{align}
    \label{eq:h_def}
h_{m,\lambda}(x) := \frac{\Po_{\lambda+x}(m)-\Po_{\lambda}(m)}{\Po_{\lambda}(m)} = 
e^{-x}(1 + \frac{x}{\lambda})^m -1. 
\end{align}
We may view $h_{m,\lambda_0}(nt)$ as the relative difference in mass probability of $X_m$ resulting from a perturbation of $U_i=1/N$ by $t$. In particular, under $H(Q)$, the mean of $X_m$ satisfies 
\begin{align}
\mu_{m}(Q) & :=  \sum_{i=1}^N \Po_{nQ_i}(m) = \mu_m^0 + \Po_{\lambda_0}(m) \sum_{i=1}^N h_{m,\lambda_0}(nQ_i-\lambda_0).
\end{align}
Likewise, under a rate prior $\pi$, 
\begin{align}
\mu_{m}(\pi) & :=  \exsub{Q\sim \pi}{\mu_m(Q)} \nonumber \\
& = \mu^0_{m} + \Po_{\lambda_0}(m) \sum_{i=1}^N \int_{\reals} h_{m,\lambda_0}(n t- \lambda_0) \pi_i(dt) \nonumber \\
& =: \mu^0_{m} + \Delta_m(\pi). \label{eq:Delta_def} 
\end{align}
Namely, $\Delta_m(\pi)$ is the expected difference in $X_m$ due to the random perturbations $\pi$ of the uniform rate sequence. Notice the identity $\sum_{m=0}^\infty \Delta_m = 0$ that follows from $\sum_{m=0}^\infty \mu_m(Q) = \sum_{m=0}^\infty \mu_m^0 = N$.

\subsection{Asymptotic normality of a family of linear tests}
\label{sec:asymptotic_normality}
The following proposition establishes conditions under which linear statistics of the form \eqref{eq:linear_test_stat} are asymptotically normal under rate sequence priors. 
\begin{prop}
\label{prop:asymp_normality}
Consider a sequence of multivariate Poisson models \eqref{eq:hyp_Q} indexed by $n$ and $N$, where $N$ and $n$ go to infinity with $N = o(n^2)$. Set 
\begin{subequations}
\label{eq:xi_r_explicit}
\begin{align}
    \xi & = \log(N)/(n \sqrt{\log(\log(N))}), \\
    r^2 & = \log(N)/(n\sqrt{N}),
\end{align}
\end{subequations}
and let $\pi \in \Pi_{\epsilon,\xi,r}$. 
Let $\{w_m\}_{m=0}^\infty$ be a non-constant sequence that satisfies the conditions:
\begin{subequations}
\label{eq:w_conditions}
\begin{equation}
    \label{eq:w_cond1}
    w_m \leq C_0e^{C_1 m },\quad \forall m\geq 0,
\end{equation}
for some $C_0$ and $C_1$ that are independent of $m$, and 
\begin{align}
&\frac{\sum_{m=0}^\infty \left|w_m\right|^4 \Po_{\lambda}(m) }{N\left( \sum_{m=0}^\infty |w_m|^2 \Po_{\lambda}(m) - \left(\sum_{m=0}^\infty w_m \Po_{\lambda}(m) \right)^2 \right)^2} \to 0.
    \label{eq:w_cond2}
\end{align}
\end{subequations}
Then the following hold. 
\begin{itemize}
\item[(i)] Uniformly in $s$,
\begin{align*}
\lim_{N \to \infty}  \Pr_{Q \sim \pi}\left[ \frac{ \sum_{m=0}^\infty w_m (X_m - \mu_m(\pi) )}{\sqrt{\langle w, \Sigma w \rangle }} \leq s\right] = \Phi\left(s \right),
\end{align*}
i.e., $T_w$ is asymptotically normal with mean $\langle w, \mu^0 + \Delta(\pi) \rangle$ and variance $\langle w, \Sigma w \rangle$. 
\item[(ii)] Set
\[
u^2(w;\pi) := \frac{\langle w,\Delta(\pi) \rangle^2 }{\langle w,\Sigma w\rangle}. 
\]
The Bayes risk of the test $\psi_{w,\alpha^*}$ of the form \eqref{eq:test_alpha_general} with threshold $z_{1-\alpha^*} = u(w;\pi)/2$ satisfies
\begin{align}
    \label{eq:alpha_opt_risk}
 \rho(\psi_{w,\alpha^*};\pi) = 2 \Phi\left(- u(w;\pi)/2 \right) + o(1),
\end{align}
and further, $\inf_{\alpha\in (0,1)} \rho(\psi_{w,\alpha};\pi) = \rho(\psi_{w,\alpha^*};\pi) + o(1)$. 
\end{itemize}
\end{prop}
The proof of Proposition~\ref{prop:asymp_normality} appears in Section~\ref{sec:proof:prop:asymp_normality}. 

We note that the explicit forms of $\xi$ and $r$ in \eqref{eq:xi_r_explicit} are chosen for convenience; they may be replaced by any sequences satisfying the general conditions \eqref{eq:max_test_cond} and \eqref{eq:L2_sep_condition}, respectively.

Condition \eqref{eq:w_cond2}, required for a Berry-Esseen-type central limit theorem, can be seen as a statement about the variations within the sequence $\{w_m\}$ when averaged under the Poisson probability mass function. This condition is invariant under rescaling of $w$, and hence is unrelated to the growth condition \eqref{eq:w_cond1}.

Under the conditions of Proposition~\ref{prop:asymp_normality}, the Bayesian hypothesis testing problem reduces, in the asymptotic limit, to a Gaussian shift experiment:
\begin{align}
\label{eq:Gaussian_shift}
\frac{T_w- \langle w,\mu^0\rangle}{\sqrt{\langle w,\Sigma w\rangle}} \overset{d}{\to} \begin{cases}
    \Ncal\left( 0, 1\right), & Q = U \\
    \Ncal\left( u(w;\pi), 1 \right), & Q \sim \pi.
\end{cases}
\end{align}
This observation motivates the construction of the following test, which is later shown to be asymptotically minimax under both the Poisson and Gaussian models.

Define the function
\begin{align}
    \label{eq:gm_def}
g_{m,\lambda}(x) & := \frac{h_{m,\lambda}(x) + h_{m,\lambda}(-x)}{2},
\end{align}
which symmetrizes $h_{m,\lambda}(x)$. Consider the test $\psi^*$ defined as in \eqref{eq:linear_test_stat}, with weights $w^* = \{w^*_m\}_{m=0}^\infty$ given by
\begin{subequations}
\label{eq:optimal_test}
\begin{align}
    \label{eq:optimal_w}
    w_m^* := g_{m,\lambda_0}(n \epsilon N^{-1/p}),\quad m=0,1,\ldots,
\end{align}
and threshold level 
\begin{align}
    \label{eq:optimal_alpha}
z_{1-\alpha} = \frac{\langle w^*,\Delta(\pi) \rangle}{2 \sqrt{\langle w^*,\Sigma w^* \rangle}}.
\end{align}
\end{subequations}
The result below characterizes the asymptotic Bayes risk of $\psi^*$ under a relevant class of priors, including a specific prior of the form \eqref{eq:three_point_prior}. 
\begin{prop}
\label{prop:risk_wstar_pi}
Consider a sequence of multivariate Poisson models \eqref{eq:hyp_Q} indexed by $n$ and $N$ with $Q\sim \pi$. Set $\xi$ and $r^2$ as in \eqref{eq:xi_r_explicit}. Assume $\epsilon N^{1-1/p} \leq \xi$. Suppose that $N$ goes to infinity and $N=o(n^2)$. 
\begin{itemize}
    \item[(i)] 
    If $\pi \in \Pi_{\epsilon,\xi,r}$, then 
    \[
\rho(\psi^*, \pi) = 2 \Phi(-u(w^*;\pi)/2) + o(1).
\]
\item[(ii)] Let $\pi^* = \prod_{i=1}^N \pi_i^*$, where 
\begin{align}
    \label{eq:pi_star_l1}
\pi_i^* = \frac{1}{2} \delta_{\frac{1}{N}+\epsilon N^{-1/p}} + \frac{1}{2} \delta_{\frac{1}{N}-\epsilon N^{-1/p}}.
\end{align}
Then $\pi^* \in \Pi_{\epsilon,\xi,r}$ for all $N$ large enough, and 
\[
\rho(\psi^*,\pi^*) = 
2 \Phi \left(- 
\sqrt{N/2} \sinh \left(n \epsilon^2 N^{1-2/p}/2\right)
\right) + o(1)
\]
\end{itemize}
\end{prop}
The proof of Proposition~\ref{prop:risk_wstar_pi} is in Section~\ref{sec:proof:prop:risk_wstar_pi}. 

Our main results, provided in the section below, say that under some conditions $\psi^*$ is asymptotically minimax and $\pi^*$ is a least favorable prior, and thus the minimax risk is asymptotically equivalent to $\rho(\psi^*,\pi^*)$. 

\section{Main Results
\label{sec:main_results}}

\subsection{Upper bound}
We obtain the upper bound by deriving the asymptotic Bayes risk under a certain test of the form \eqref{eq:linear_test_stat} with weights provided in \eqref{eq:optimal_w}. Maximizing this risk over alternative priors results in $\pi^*$ of \eqref{eq:pi_star_l1} as the unique maximizer. Consequently, the asymptotic Bayes risk of $\pi^*$ bounds the minimax risk from above. 

\begin{theorem}
\label{thm:upper_bound}
Consider the multivariate Poisson model \eqref{eq:hyp_Q} under the minimax setting \eqref{eq:minimax_risk}. As $N\to \infty$ with $N=o(n^2)$,
\[
R^*(V_\epsilon) + o(1) \leq 
2\Phi\left(- 
\sqrt{\frac{N}{2}} \sinh \left( \frac{n \epsilon^2 N^{1-2/p}}{2}\right)
\right). 
\]
\end{theorem}
The proof of Theorem~\ref{thm:upper_bound} is in Section~\ref{sec:proof:upper_bound}. 

\subsection{Lower bound}
We now use the prior $\pi^*$  of \eqref{eq:pi_star_l1} to bound the minimax risk from below. This is achieved by analyzing the power of the likelihood ratio test for testing $H(U)$ against the simple alternative $H_1 \,:\, Q \sim \pi^*$. 

\begin{theorem}
\label{thm:lower_bound}
Consider the multivariate Poisson model \eqref{eq:hyp_Q} under the minimax setting \eqref{eq:minimax_risk} where $N$ and $n$ tend to infinity. If $\log(N) / (n N^{1-1/p} \epsilon) = O(1)$ then $R^*(V_\epsilon)=o(1)$. Otherwise, assuming $N=o(n^2)$, 
\begin{align}
    \label{eq:prf:lower_bound}
        2 \Phi(- \tilde{u}/2) \leq R^*(V_{\epsilon}) + o(1),
\end{align}    
where 
\begin{align*}
    \tilde{u}^2 := \frac{\langle  \tilde{w}, \Delta(\pi^*) \rangle^2}{\langle \tilde{w}, \Sigma \tilde{w}\rangle}.
\end{align*} 
\end{theorem}
The proof of Theorem~\ref{thm:lower_bound} is provided in Section~\ref{sec:proof:lower_bound}.

\begin{figure*}[t]
    \centering
\includegraphics[scale=.45]{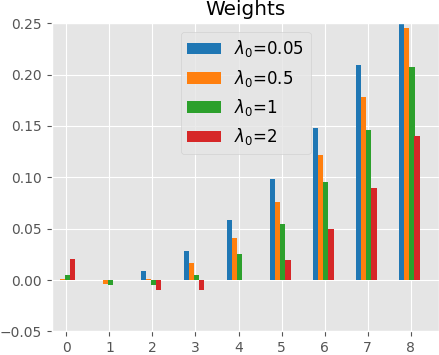}
~
\includegraphics[scale=.45]{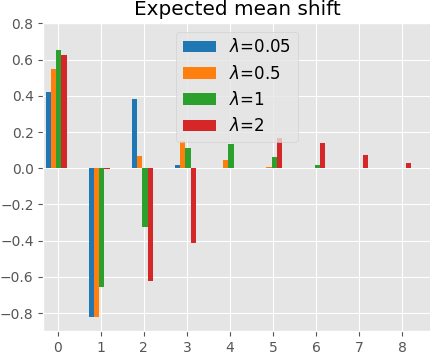}
~
\includegraphics[scale=.45]{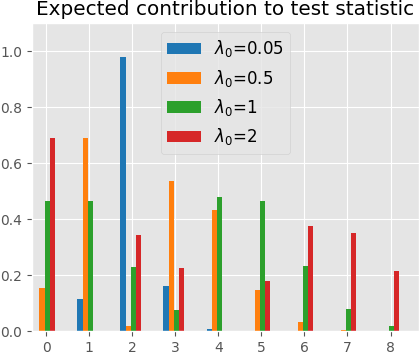}
    \caption{
    The weights of the minimax test $w^*$ of \eqref{eq:optimal_test} and their relative expected contributions to departures from the null in the test statistic; here $n=10,000$, $p=1$, $\epsilon=0.1$. Left: Bars proportional to the coordinates of $w^*$, where $w^*_0$ is the weight for missing categories, $w^*_1$ to singletons, $w^*_2$ to exclusive collisions, and so on. Center: The expected shift in the mean of the histogram ordinates under the least favorable prior $\pi^*$ defined in \eqref{eq:pi_star_l1}.
    Right: the normalized product of $w_m^* \Delta_m(\pi^*)$ indicating the expected difference each histogram ordinate contributes to the minimax test statistic under $\pi^*$ relative to the null. Different colors represent different values of $\lambda_0 = n/N$. As $\lambda_0 \to 0$, only exclusive collision count $X_2$ contributes to the test statistic, as indicated by the blue bar in the right panel. 
    }
    \label{fig:kernel}
\end{figure*}

\subsection{The minimax risk}
The upper bound from Theorem~\ref{thm:upper_bound} and the lower bound from Theorem~\ref{thm:lower_bound} coincide asymptotically whenever the minimax risk does not vanish. 

\begin{theorem}
\label{thm:matching}
Consider the multivariate Poisson model \eqref{eq:hyp_Q} under the minimax setting \eqref{eq:minimax_risk} with $N$ and $n$ tend to infinity. Suppose that $N = o(n^2)$ and 
\begin{align}
    \label{eq:thm:matching_cond}
    \lim_{N \to \infty} R^*(V_\epsilon) \in (0,1].
\end{align}
Then 
\begin{align}
    \lim_{N \to \infty} \frac{R^*(V_\epsilon)}{2 \Phi(-u_{\epsilon,n,N,p}/2)} = 1.
    \label{eq:thm:matching}
\end{align}
\end{theorem}
The proof of Theorem~\ref{thm:matching} is in Section~\ref{sec:proof:matching}. Sufficient conditions for \eqref{eq:thm:matching_cond} follow from Theorem~\ref{thm:lower_bound} and the observation that if $\lim_{N\to \infty} u_{\epsilon,n,N,p} > 0$, then $\epsilon N^{1-1/p}=o(1)$ if and only if $N=o(n^2)$. 
\begin{cor}
\label{cor:matching}
Suppose that $u_{\epsilon,n,N,p} \to c$ for some $c \in (0,\infty)$. If $\epsilon N^{1-1/p} = o(1)$ or $N=o(n^2)$, then 
\begin{align}
     R^*(V_\epsilon) = 2 \Phi(-c/2) + o(1).
    \label{eq:cor:matching}
\end{align}
\end{cor}
When $u_{\epsilon,n,N,p} \to \infty$, Theorem~\ref{thm:upper_bound} implies that $R^*(V_\epsilon) = o(1)$, so both $R^*(V_\epsilon)$ and $2\Phi(-u_{\epsilon,n,N,p}/2)$ vanish. However, in this regime $R^*(V_\epsilon)$ may decay faster, and thus \eqref{eq:thm:matching} need not hold. 

We note that when $u_{\epsilon,n,N,p}\to \infty$, Mill's ratio implies (c.f. \cite{shorack2009empirical}) 
 \begin{align*}
2\Phi\left( -\frac{u_{\epsilon,n,N,p}}{2}\right) = \frac{2\exp\{\frac{-(u_{\epsilon})^2}{8}\} (1+o(1))}{u_{\epsilon,n,N,p} \sqrt{ \pi/2}},
\end{align*}
which leads to 
\begin{align}
\label{eq:sample_complexity_uniform}
    \frac{4\sqrt{N \log(1/R^*(V_{\epsilon}))}}{n \epsilon^2} = 1 + o(1).
\end{align}
The scaling of $n$ deduced from \eqref{eq:sample_complexity_uniform} is equivalent to the one suggested in \cite{gupta2022sharp}. However, the approximation to $R^*(V_{\epsilon})$ obtained by the non-vanishing terms in \eqref{eq:sample_complexity_uniform} is loose unless $R^*(V_{\epsilon})$ is very small.

\subsection{The minimax test}
\label{sec:minimax_test}
It follows from Theorem~\ref{thm:matching} that when $\lim_{N \to \infty} R^*(V_\epsilon) > 0$, an asymptotically least favorable prior is given by $\pi^*$ of \eqref{eq:pi_star_l1}, as conceptually illustrated by the 4 red dots in Figure~\ref{fig:alternative}. 
Additionally, the test $\psi^*$ of \eqref{eq:optimal_test} is asymptotically minimax under the assumptions of Theorem~\ref{thm:matching}.

The left panel in Figure~\ref{fig:kernel} illustrates the weights $\{w^*_m\}$ for several values of $\lambda_0$ and $m$ for $p=1$. By a second-order approximation to $g_{m,\lambda_0}$ (see Lemma~\ref{lem:gm_properties}), these weights satisfy
\begin{align*}
w^*_m =  \frac{1}{2}(\epsilon N^{1-1/p})^2\left( \left( m - \lambda_0 \right)^2 - m + o(1) \right).
\end{align*}
This shows that the minimax test behaves similarly to the chi-squared test for large values of $m$ due to the quadratic term, but with some differences for small $m$ that are particularly apparent when $\lambda_0$ approaches zero. To understand these differences, we analyze below the asymptotic risk of the chi-squared and collision-based tests under the least favorable prior $\pi^*$. 

The expected difference between the null and alternative due to the $m$-th histogram ordinate is given by 
\[
T^*_m := w_m^*\Delta_m(\pi^*) = N \Po_{\lambda_0}(m)\left(g_{m,\lambda_0}(n \epsilon N^{-1/p}) \right)^2.
\]
For small $\lambda_0 = n/N$ with $ \lambda_0 = o(N^{1-1/p} \epsilon)$, we have $T_2 \propto \lambda_0^2(1+o(1))$ while $T_m \propto o(\lambda_0^2)$ for $m \neq 2$. It follows that in this limit the minimax test statistic is only affected by $T^*_2$, so the exclusive collision statistic $X_2$ accounts for all the difference. The situation is illustrated by the panel on the right-hand side in Figure~\ref{fig:kernel}. 

\subsection{The asymptotic risk of chi-squared and collision statistics}
The chi-squared test statistic defined in \eqref{eq:chisquared_def} satisfies
\begin{align}
     T_{\chi^2}= \sum_{m=0}^\infty \frac{(m-\lambda_0)^2}{\lambda_0}X_m, 
\end{align}
so it is of the form \eqref{eq:linear_test_stat} with weights
\begin{align}
    \label{eq:w_chisquared}
w_{\chi^2} := \left(\lambda_0, \frac{(1-\lambda_0)^2}{\lambda_0},\frac{(2-\lambda_0)^2}{\lambda_0},\ldots \right).
\end{align}
Similarly, the collision statistic is given by 
\begin{align}
    \label{eq:collision_test}
    T_{\col} := \sum_{i=1}^N \binom{O_i}{2} = \sum_{m=0}^{\infty} X_m w_m,
\end{align}
with $w_m = m(m-1)/2$. Under the least favorable prior $\pi^*$ given in \eqref{eq:pi_star_l1}, the asymptotic power of the test $\psi_{\chi^2}$, based on $T_{\chi^2}$, and the test $\psi_{\col}$, based on $T_{\col}$, can be derived using arguments similar to those in the proof of Theorem~\ref{thm:upper_bound}, as stated below.
\begin{prop}
\label{prop:chisq_n_collision}
    As $N \to \infty$, suppose that $u_{\epsilon,n,N,p} \to c$ for some $c \in (0,\infty)$ and $\epsilon N^{1-1/p} \to 0$ or $N=o(n^2)$. Then
    \begin{align}
        \label{eq:asymp_risk_chisq}
        \rho(\psi_{\chi^2},\pi^*) = 2 \Phi \left(-\sqrt{\frac{\lambda_0}{1+\lambda_0}}\frac{u_{\epsilon,n,N,p}}{2}\right)+ o(1),
    \end{align}
    and 
    \begin{align}
        \label{eq:asymp_risk_col}
        \rho(\psi_{\col},\pi^*) = 2 \Phi \left(-\sqrt{\frac{1}{1+2\lambda_0}}\frac{u_{\epsilon,n,N,p}}{2}\right)+ o(1),
    \end{align}
\end{prop}
The proof of Proposition~\ref{prop:chisq_n_collision} appears in Section~\ref{sec:chisq_n_collision:proof}.

Collision-based tests have been proposed and analyzed in several works in settings related to ours, particularly when $\lambda_0 \ll 1$ \cite{goldreich2000testing,batu2001testing,rubinfeld2005testing,diakonikolas2019collision}. While tests based solely on $X_2$ (exclusive collisions) are sometimes considered, they are generally less powerful than the full collision statistic in \eqref{eq:collision_test}; see \cite{diakonikolas2019collision}. 



\begin{figure}
    \centering
    \includegraphics[scale=1]{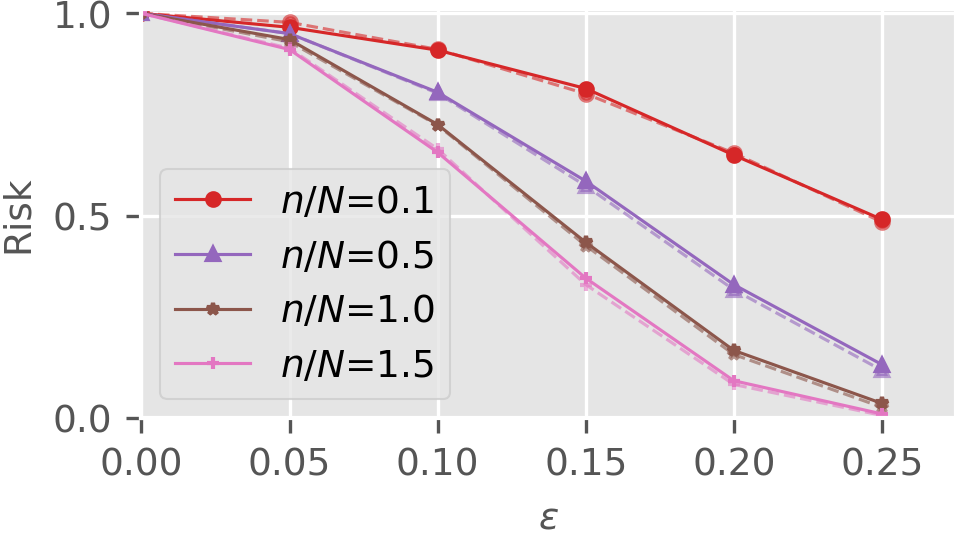}
    \caption{
        Empirical (continuous) and theoretical (dashed) risk under the least favorable prior 
        versus $\epsilon$ for several values of $\lambda_0 = n/N$ and $p=1$. The empirical risk in each configuration is the average error in $10,000$ Monte-Carlo trials. In each trial, we used $n=10,000$ samples from the null and $n=10,000$ samples from the alternative to evaluate the Type-I and Type-II errors, respectively. }
        \label{fig:risk}
    \end{figure}

\section{Empirical Comparisons}
\label{sec:numerical_evalutations}

We compare the asymptotic minimax risk $R^*(V_{\epsilon})$, obtained from \eqref{eq:thm:matching}, to an empirically estimated risk computed via Monte-Carlo simulations under the least favorable prior $\pi^*$ from \eqref{eq:pi_star_l1} and under the null in Figure~\ref{fig:risk}. For each configuration, we estimate the empirical risk as the sum of the Type I error rate (the proportion of rejections in $\nMonte$ independent trials under the null) and the Type II error rate (the proportion of non-rejections in $\nMonte$ independent trials under the alternative). We also evaluate the empirical risks of the chi-squared test from \eqref{eq:chisquared_def} and of a collision-based test defined in \eqref{eq:collision_test}.

\begin{figure}
    \label{fig:comparisons}
    \centering
    \includegraphics[scale=1]{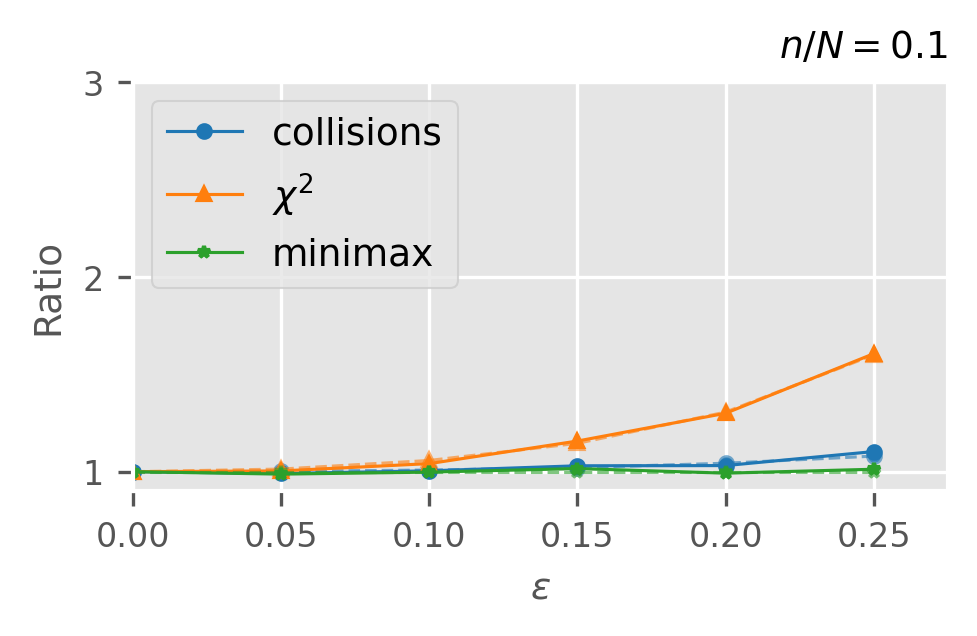}

    \includegraphics[scale=1]{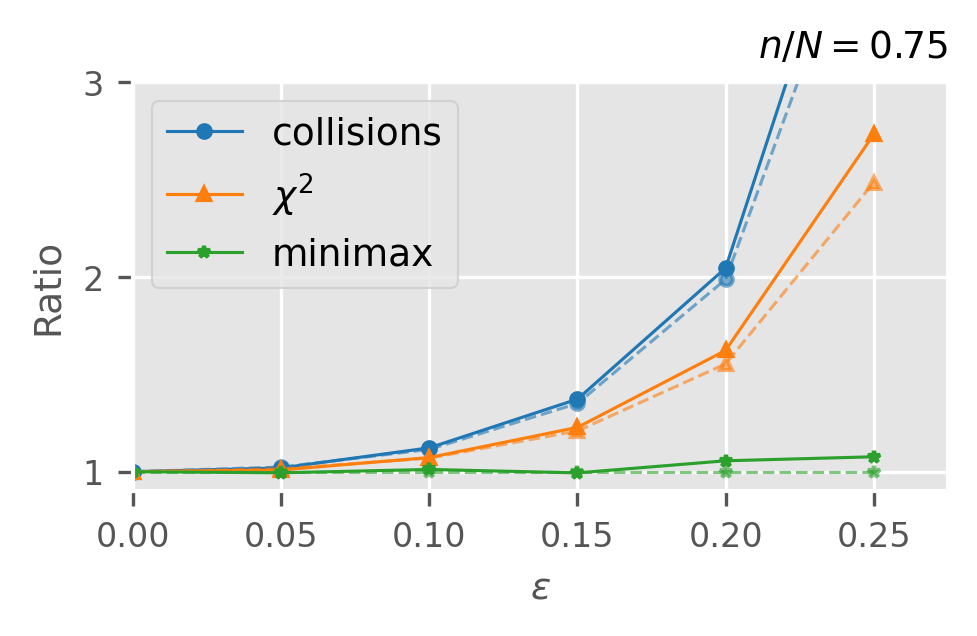}
\caption{Risk of the chi-squared test \eqref{eq:chisquared_def} and the collision-based test \eqref{eq:collision_test}, normalized by the asymptotic minimax risk $R^*(V_{\epsilon})$, under the least favorable prior and $p=1$. Each curve shows the ratio to $R^(V_{\epsilon})$; solid lines indicate empirical risk, and dashed lines indicate theoretical asymptotic risk. Top: $\lambda_0 = n/N = 1/10$. Bottom: $\lambda_0 = 3/4$. 
Empirical risks are computed as the average error over $\nMonte$ Monte Carlo trials. In each trial, $n=\nval$ samples are drawn from the null and $n=\nval$ samples from the alternative to estimate the Type-I and Type-II errors, respectively.
}
\end{figure}

In Fig.~\ref{fig:comparisons}, we report the empirical risks of the chi-squared and collision-based tests under the least favorable prior, together with their theoretical asymptotic risks derived in Proposition~\ref{prop:chisq_n_collision}. The results show that the minimax test uniformly outperforms both alternatives across all configurations considered. In addition, when $\epsilon$ is small, the asymptotic risk of each test provides an accurate approximation to its empirical risk.

\begin{figure}
\label{fig:constant_risk}
    \centering
    \includegraphics[scale=1]{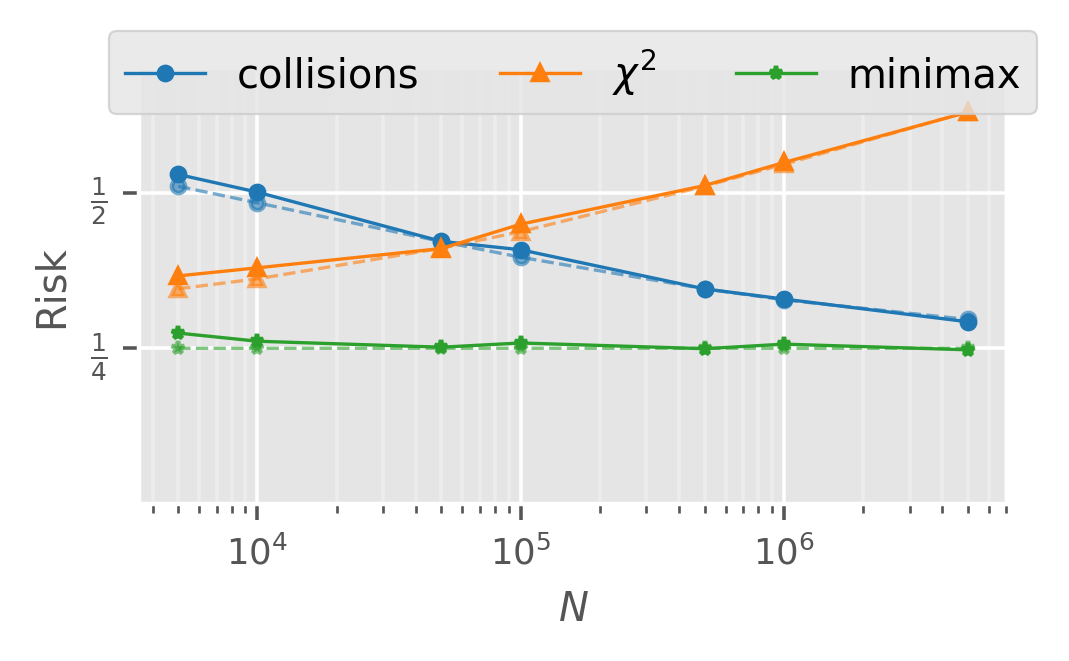}

    \includegraphics[scale=1]{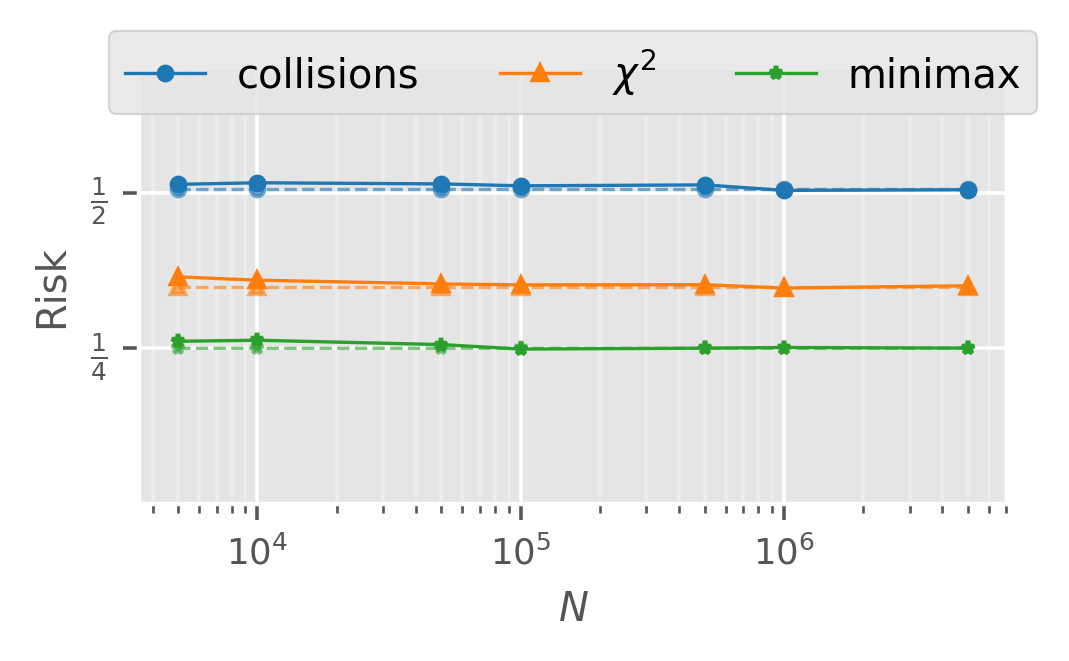}
\caption{Empirical risk of the minimax test under the least favorable prior, together with the risks of the chi-squared and collision-based tests, with parameters scaled in $N$ to yield a constant minimax risk; dashed lines indicate the theoretical asymptotic risks. Top: $\lambda_0 \to 0$ as $N \to \infty$. Bottom: $\lambda_0 = 1$. In both panels, $p=1$. Theoretical risks are evaluated according to Theorem~\ref{thm:matching} and Proposition~\ref{prop:chisq_n_collision}; empirical risks are computed as in Fig.~\ref{fig:comparisons}.
}
\end{figure}

A similar comparison is presented in Fig.~\ref{fig:constant_risk}, where all parameters are scaled so as to attain a constant minimax risk. When $\lambda_0 \to 0$, the risk of the collision-based test converges to the minimax risk. In contrast, when $\lambda_0$ is fixed, the asymptotic risks of both the minimax and chi-squared tests under the least favorable prior remain separated from the minimax risk by a constant gap. These behaviors are consistent with the asymptotic risk expressions in Proposition~\ref{prop:chisq_n_collision}.



\section{Discussions}
\label{sec:discussion}

\subsection{Multinomial sampling}

As noted in the introduction, the minimax risk under the Poisson model serves as an upper bound for the corresponding risk under multinomial sampling. Consequently, the upper bound derived in Theorem~\ref{thm:upper_bound} applies directly to the multinomial setting. In a separate note \cite{kipnis2026CLT}, we show that the lower bound in Theorem~\ref{thm:lower_bound} also holds asymptotically under multinomial sampling. Thus, the asymptotic minimax risk in either model is fully characterized by Theorem~\ref{thm:matching}. 

The primary challenge in extending the lower bound to the multinomial model lies in establishing a conditional central limit theorem for the test statistic under the constraint $\sum_{i=1}^N O_i = n$ and arguing that its variance is unchanged. While this de-Poissonization step is classical \cite{renyi1962three,holst1972asymptotic,morris1975central,quine1984normal,ingster1993asymptotically,Nussbaum1996,ermakov1998asymptotic,janson2001moment,barbour2009small}, the high-dimensional regime and the specific alternatives considered here require a refinement of existing results. In the companion note \cite{kipnis2026CLT}, we extend the techniques of \cite{janson2001moment,klein2015conditional} to derive the necessary conditional limit theorem for the likelihood ratio statistic associated with the least favorable prior $\pi^*$ in \eqref{eq:pi_star_l1}. In particular, the (unconditioned) minimax statistic $T_w^*$ is asymptotically minimax under the multinomial model in the primary regime of interest $u_{\epsilon,n,N,p} = O(1)$. In fact, in this setting, $T_{w^*}$, the chi-squared test statistic $T_{\chi^2}$, and the collision statistic $T_{\col}$ are all asymptotically minimax. However, out of the three, only the distribution of $T_{w^*}$ is unchanged under the sample size conditioning, which is key for establishing the asymptotic minimax risk under multinomial sampling. 


\subsection{Removal of non-ball shapes}
The similarities between some parts of our analysis to those of \cite{ingster1987minimax,suslina1999extremum} suggest the possibility of generalizing our setting to the removal of other geometric shapes like ellipsoids and Besov bodies  \cite{ingster2000minimax}. Such extensions may lead to multi-dimensional versions of the optimization problem considered  
in the proof of Theorem \eqref{thm:upper_bound} in which the parameter is a pair of positive sequences $\{\eta_i\}_{i=1}^N$ and $\{\mu_i\}_{i=1}^N$ rather than a pair of positive numbers. 

\subsection{Sparse alternatives}
Another interesting extension is attained by imposing on $V_\epsilon$ constraints of the form $\|Q_i - U_i\|_q \leq \epsilon'$ for some $q>0$ and $\epsilon'>0$. For example, the case of $q$ close to zero is related to sparse alternatives as considered in \cite{bhattacharya2021sparse} and \cite{arias2015sparse}. The methods developed in this paper appears suitable for analyzing such extensions.


\section{Proofs\label{sec:proofs}}
This section contains the proofs of 
Lemma~\ref{lem:R_geq_rho}, 
Proposition~\ref{prop:asymp_normality}, Proposition~\ref{prop:risk_wstar_pi}, 
Theorem~\ref{thm:upper_bound}, Theorem~\ref{thm:lower_bound}, Theorem~\ref{thm:matching}, and Proposition~\ref{prop:chisq_n_collision}

\subsection{Proof of Lemma~\ref{lem:R_geq_rho}
\label{sec:R_geq_rho:proof}}

Because every sequence $Q \in V_{\epsilon,\xi,r}$ corresponds to a prior whose $i$-th coordinate is a mass probability $Q_i$, we have $\rho^*(\Pi_{\epsilon,\xi,r}) \geq R^*(V_{\epsilon,\xi,r})$. It is left to show the reversed inequality. Fix $\pi \in \Pi_{\epsilon,\xi,r}$. Denote by $P_\pi$ the probability law $Q \sim \pi$. Each $\pi_i$ is supported in $[1/N-\xi,1/N+\xi]$, hence the variance of $(Q_i-1/N)^p$ exists and uniformly bounded in $i$. 

Denote $V = V_1 \times ... \times V_N$. Because $\pi = \prod_{i=1}^N \pi_i$, the law of large number implies
\begin{align*}
\liminf_{N \to \infty} & \frac{1}{N}\sum_{i=1}^N \epsilon^{-p}\left|NQ_i - 1\right|^p \\
&\qquad = \liminf_{N \to \infty} \frac{\epsilon^{-p}}{N}\sum_{i=1}^N \ex{\left|NQ_i - 1\right|^p} \geq 1
\end{align*}
$P_\pi$ almost surely. Likewise,
\begin{align*}
\limsup_{N \to \infty} & \frac{1}{N}\sum_{i=1}^N r^{-2}N\left|Q_i - 1/N\right|^2 \\
& \overset{a.s.}{=} \liminf_{N \to \infty} \frac{r^{-2}}{N}\sum_{i=1}^N N \ex{\left|Q_i - 1/N\right|^2} \leq 1. 
\end{align*}
If follows that for some sequences $r_N$ and $\epsilon_N$ satisfying
$\epsilon_N/\epsilon \to 1$, $r_N/r \to 1$, 
$\epsilon_N/\epsilon \geq 1$, and $r_N/r \geq 1$,  
$P_{\pi}\left[Q \in V_{\epsilon_N,\xi,r_N}\right] \to 1$. We now argue similarly to \cite[Ch. 4.1]{ingster1993asymptotically}: Define the measure $\tilde{\pi}(A) = \pi(A \cap V_{\epsilon_N,\xi,r_N})/\pi(V_{\epsilon_N,\xi,r_N})$. We have $\tilde{\pi}(V_{\epsilon_N,\xi,r_N})=1$, and hence $\rho^*(\tilde{\pi}) \leq R^*(V_{\epsilon_N,\xi,r_N})$. On the other hand, the total variation distance obeys $\mathrm{TV}(P_\pi,P_{\tilde{\pi}})\to 0$ 
because $\pi(V_{\epsilon_N,\xi,r_N}) \to 1$, implying $\rho^*(\pi) = \rho^*(\tilde{\pi})+o(1) \leq R^*(V_{\epsilon_N,\xi,r_N})+o(1) = R^*(V_{\epsilon,\xi,r})+o(1)$, the last transition by continuity of $R^*(V_{\epsilon,\xi,r})$ in $\epsilon$ and $r$. 
\qed

\subsection{\label{sec:proof:prop:asymp_normality} Proof of Proposition~\ref{prop:asymp_normality}}
We first state and prove a similar result for rate sequences $Q$ that are sufficiently close to $U$.
\begin{lem}
\label{lem:asym_normality1}
Consider a multivariate Poisson model
\[
O_i \sim \Pois(n Q_i)\, \text{ independently for } i=1,\ldots,N, 
\]
where $N$ and $n$ go to infinity. Let $U =(\lambda,...,\lambda) \in \reals^N$, $\lambda=n/N$. Suppose that $Q \in \reals_+^N$ satisfies:
\begin{subequations}
\label{eq:Q_cond}
\begin{align}
    \left\| Q - U \right\|_2^2 \leq a_{N} \quad \text{and} \quad
     \max\{ \max_{i} Q_i, 1/N\} & \leq b_{N}, 
\end{align}
for some sequences $\{a_{N}\}$ and $\{b_{N}\}$ such that
\begin{align}
\label{eq:an_bn_cond}
n^2 N^{-1} a_{N} e^{C n b_{N}} =o(1),
\end{align}
for any fixed $C>0$. 
\end{subequations}
Let $\{w_m\}_{m=0}^\infty$ be a non-constant sequence that satisfy \eqref{eq:w_conditions}. 
Then, uniformly in $s$,  
\begin{align}
\lim_{n \to \infty}
\Prp{\frac{ \sum_{m=0}^\infty w_m (X_m - \mu_m(Q))}{ \sqrt{ \langle w,\Sigma w \rangle  } } \leq s } = \Phi\left(s \right)
\label{eq:lem1:statement}
\end{align} 
\end{lem}
The proof of Lemma~\ref{lem:asym_normality1} is in the Appendix (\ref{sec:proof:lem:asym_normality1}). 

We use Lemma~\ref{lem:asym_normality1} with sequences 
\begin{align*}
b_N& =\frac{1}{N} + \xi = \frac{1}{N} + \frac{\log(N)}{n \sqrt{\log(\log(N))}}, \\
a_N & = r^2(1+\delta_N) = (1+\delta_N)\frac{\log(N)}{n \sqrt{N}},
\end{align*}
for some $\delta_N \to 0$ that we specify later. Notice that these $b_N$ and $a_N$ satisfy 
\eqref{eq:an_bn_cond}. Lemma~\ref{lem:asym_normality1} applies, conditioned on the event 
\[
A_{N} := \{ \max_{Q_i} \leq b_N~~  \text{and}~~
\|Q -U\|_2^2 \leq a_N\}.
\]
Namely, 
\begin{align*}
\lim_{N \to \infty}
\Pr_{Q \sim \pi} \left[
\frac{ \sum_{m=0}^\infty w_m (X_m - \mu_m(Q))}{ \sqrt{ \langle w,\Sigma w \rangle  } } \leq s \mid A_N \right] = \Phi\left(s \right). 
\end{align*} 
It is left to show that $\Prp{A_N} \to 1$. For $Q  \sim \pi \in \Pi_{\epsilon,\xi,r}$, we have $\Prp{Q_i\leq b_N}=1$ because $\pi_i([1/N-\xi,1/N+\xi])=1$. Additionally, by the law of large numbers, almost surely,
\begin{align*}
\limsup_{N \to \infty} & \frac{1}{N} \sum_{i=1}^N r^{-2} N \left|Q_i - 1/N \right|^2 \\
& = \limsup_{N \to \infty} \frac{r^{-2}}{N} \sum_{i=1}^N N \ex{ \left|Q_i - 1/N \right|^2} \leq 1.
\end{align*}
It follows that there exists $\delta_N \to 0$ such that $\Prp{\|Q -U\|_2^2 \leq a_N}\to 1$, hence $\Prp{A_N} \to 1$. Part (i) follows. 

By (i), we conclude that the power of the test \eqref{eq:test_alpha_general} converges to $
\Phi\left(-(u(w;\pi) - z_{1-\alpha}) \right)$, where
\[
u^2(w;\pi) := \frac{\langle w,\Delta(\pi) \rangle^2 }{\langle w,\Sigma w\rangle}. 
\]
The Bayes risk of $\psi_{w,\alpha}$ satisfies 
\begin{align*}
o(1) + \rho(\psi_{w,\alpha}) & = \alpha + \Phi\left(-(u(w;\pi) - z_{1-\alpha}) \right) \\
& = \Phi(-z_{1-\alpha}) + \Phi\left(-(u(w;\pi) - z_{1-\alpha}) \right).
\end{align*}
The last expression is minimal when $z_{1-\alpha} = u(w;\pi)/2$. This implies (ii). 
\qed

\subsection*{Proof of Proposition~\ref{prop:risk_wstar_pi}
\label{sec:proof:prop:risk_wstar_pi}}
Below are some useful properties of $g_{m,\lambda}(x)$ that follow from Lemmas~\ref{lem:gm_properties} and \ref{lem:gm_sums} in the Appendix.
\begin{itemize}
    \item[(i)] $\sum_{m=0}^\infty g_{m,\lambda}(t) \Po_{\lambda}(m) = 0$ for any $t$ 
    \item[(ii)] $\sum_{m=0}^\infty g_{m,\lambda}^2(t) \Po_{\lambda}(m)=2 \sinh^2(t^2/(2\lambda))$
    \item[(iii)] $\sum_{m=0}^\infty g_{m,\lambda}^4(t) \Po_{\lambda}(m) = F(t^2/\lambda)$ for some non-negative and continuous function $F$. Furthermore, for $t = o(\sqrt{\lambda})$
    \[
    F\left( \frac{t^2}{\lambda}\right) =  \frac{(t^2/\lambda)^4}{2\lambda^2}(1+o(1)).
    \]
\end{itemize}
Set $\lambda_0 = \lambda = n/N$,  $t_N := t = n \epsilon N^{-1/p}$ and $w_m = g_{m,\lambda_0}(t_N)$. The assumption $\epsilon N^{1-1/p} \leq \xi$ implies $(n \epsilon N^{-1/p})^2/\lambda_0 = o(1)$. By properties (i)-(iii) above,
\begin{align}
    & \frac{\sum_{m=0}^\infty \left|w_m\right|^4 \Po_{\lambda_0}(m)  }{N\left( \sum_{m=0}^\infty |w_m|^2 \Po_{\lambda_0}(m) - \left(\sum_{m=0}^\infty w_m \Po_{\lambda_0}(m) \right)^2 \right)^2} \nonumber \\
    & \qquad = \frac{F(t_N^2/\lambda_0)}{N\left(2\sinh^2\left(\frac{t_N^2}{2\lambda_0}\right) -0\right)^2}  \label{eq:En_gm}
\end{align}
If $\lambda_0$ is bounded away from zero, $F(t_N^2/\lambda_0)$ is bounded from above, $\sinh^2\left(\frac{t_N^2}{2\lambda_0}\right)$ is bounded from below, and thus \eqref{eq:w_cond2} holds. If $\lambda_0 \to 0$, we get
\begin{align*}
\frac{F(t_N^2/\lambda_0)}{N\left(2\sinh^2\left(\frac{t_N^2}{2\lambda_0}\right)\right)^2} &= \frac{t_N^8}{2\lambda_0^6} \frac{1+o(1)}{4N \left(\frac{t_N^2}{2\lambda_0}\right)^4(1+o(1))} \\
&= \frac{2N(1+o(1))}{n^2} = o(1),
\end{align*}
by the assumption $N=o(n^2)$. Therefore, \eqref{eq:w_cond2} holds. (i) follows immediately from Proposition~\ref{prop:asymp_normality}-(ii). 
For (ii), notice that
\[
\Delta_m(\pi^*) = N\Po_{\lambda_0}(m) g_{m,\lambda_0}(n\epsilon N^{-1/p}), 
\]
hence 
\begin{align}
    \label{eq:wstar_delta_star}
\langle w^*, \Delta(\pi^*) \rangle = N \sum_{m=0}^\infty \Po_{\lambda_0}(m) g_{m,\lambda_0}^2(n\epsilon N^{-1/p}).
\end{align}
A closed-form expression to sum in \eqref{eq:wstar_delta_star} follows from Lemma~\ref{lem:gm_properties}, leading to 
\begin{align}
    \sum_{m=0}^\infty \Po_{\lambda_0}(m) g_{m,\lambda_0}^2(n\epsilon N^{-1/p}) =  2\sinh^2\left(n \epsilon^2 N^{1-2/p}/2\right).
    \label{eq:sinh}
\end{align}
Additionally, 
\[
\langle w^*, \mu^0 \rangle = \sum_{m=0}^\infty g_{m,\lambda_0}( n\epsilon N^{-1/p} ) \Po_{\lambda_0}(m) = 0, 
\]
also by Lemma~\ref{lem:gm_properties}. It follows from \eqref{eq:Sigma_matrix} that 
\begin{align}
    \label{eq:wstar_sigma_wstar}
\langle w^*, \Sigma w^* \rangle = N \sum_{m=0}^\infty \Po_{\lambda_0}(m) g_{m,\lambda_0}^2(n\epsilon N^{-1/p}) = \langle w^*, \Delta(\pi^*) \rangle.
\end{align}

Suppose first that $\pi^* \in \Pi_{\epsilon,\xi,r}$. By (i), 
\begin{align*}
    \rho(\psi^*,\pi^*) = 2 \Phi(-u(w^*;\pi^*)/2) + o(1),
\end{align*}
where \eqref{eq:wstar_delta_star}, \eqref{eq:wstar_sigma_wstar}, and \eqref{eq:sinh} imply
\begin{align*}
u^2(w^*;\pi^*) & = \frac{ \langle w^*, \Delta(\pi^*) \rangle^2}{\langle w^*, \Sigma w^* \rangle} = \langle w^*, \Sigma w^* \rangle \\
& = 2 N \sinh^2 \left(n \epsilon^2 N^{1-2/p}/2\right).
\end{align*}

It is left to show $\pi^* \in \Pi_{\epsilon,\xi,r}$. The assumption 
$\epsilon N^{1-1/p} \leq \xi$ implies 
$\epsilon N^{-1/p} - 1/N \leq \xi$ for all $N$ large enough, thus $\pi^*_i([1/N-\xi, 1/N + \xi])=1$. Additionally,
\begin{align*}
\sum_{i=1}^N \exsub{Q_i \sim \pi_i}{\left|Q_i - 1/N \right|^p} &= N \frac{\left|N^{-1/p}\epsilon\right|^p + \left|N^{-1/p} \epsilon\right|^p}{2} = \epsilon^p,
\end{align*}
and, also by $\epsilon N^{1-1/p} \leq \xi$, 
\begin{align*}
& \sum_{i=1}^N \exsub{Q_i \sim \pi_i}{\left|Q_i - 1/N \right|^2} = N^{1-2/p} \epsilon^2 \\
& \quad \leq \xi^2/N \leq \frac{\log^2(N)}{n^2 N \log(\log(N))} \\
& \quad = \frac{r^2 \log(N)}{n \sqrt{N} \log(\log(N))} \leq r^2. 
\end{align*}
for all $N$ large enough.  \qed

\subsection*{Proof of Theorem~\ref{thm:upper_bound}}
\label{sec:proof:upper_bound}
Consider the set $V_{\epsilon,\xi,r}$ of \eqref{eq:V_eps_xi_r_def} with $\xi = \log(N)/(n \sqrt{\log(\log(N))})$ and $r^2 = \log(N)/(n\sqrt{N})$. By Corollary~\ref{cor:trivial_separation}-(iii), it is enough to show 
\begin{align}
\label{eq:prf:upper_bound:to_show}
R^*(V_{\epsilon,\xi,r}) \leq 2\Phi\left(- 
\sqrt{N/2} \sinh \left(n \epsilon^2 N^{1-2/p}/2\right)
\right) + o(1)
\end{align}

If $V_{\epsilon,\xi,r}$ is empty then $R^*(V_{\epsilon,\xi,r})=0$ by convention and \eqref{eq:prf:upper_bound:to_show} holds. We henceforth assume that $V_{\epsilon,\xi,r}$ is non-empty, and thus we have \eqref{eq:non_empty_constraint} and in particular $\epsilon N^{1-1/p} \leq \xi$. Therefore, by Proposition~\ref{prop:risk_wstar_pi}-(i), $\pi \in \Pi_{\epsilon,\xi,r}$ and 
\begin{align*}
    \rho(\psi^*;\pi) = 2 \Phi\left( - u(w^*;\pi)/2 \right)+o(1),
\end{align*} 
where 
\begin{align*}
    u^2(w^*;\pi) & := \frac{ \langle w^*,  \Delta(\pi) \rangle^2}{\langle w^* , \Sigma w^* \rangle}. 
\end{align*}
We now consider the minimization of $u^2(w^*;\pi)$ over $\pi \in \Pi_{\epsilon,\xi,r}$. Notice that $\Delta_m(\pi)$ is linear in $\pi$ and  $\Pi_{\epsilon,\xi,r}$ 
is a convex set, hence the problem of minimizing $u^2(w^*;\pi)$ over $\pi \in \Pi_{\epsilon,\xi,r}$ is a convex optimization problem. Minimizing $u^2(w^*;\pi)$ over the larger set
\begin{align*}
    \Pi_{\epsilon,\xi}  := \Biggl\{ \pi = \prod_{i=1}^N \pi_i \,:& \, \sum_{i=1}^N \exsub{Q_i \sim \pi_i}{\left|Q_i-1/N \right|^{p}} \geq \epsilon^p, \nonumber \\
    & \pi_i([1/N-\xi, 1/N+\xi])=1 \Biggr\}
\end{align*}
is somewhat simpler and leads to an equivalent result because the minimizer turns out to be in the smaller set $\Pi_{\epsilon,\xi,r}$. 
The solution to this minimization is provided in the following lemma. 
\begin{lem}
\label{lem:minimizing_u}
Define the product prior $\pi^* := \prod_{i=1}^N \pi^*_i $, where 
\begin{align}
    \label{eq:pi_star_proof}
\pi_i^* =  
\frac{1}{2} \delta_{\frac{1}{N}+\epsilon N^{-1/p}} + \frac{1}{2} \delta_{\frac{1}{N}-\epsilon N^{-1/p} },\quad i=1,\ldots,N.
\end{align}
Then
\begin{align*}
 \inf_{\pi \in \Pi_{\epsilon,\xi}} \langle w^*, \Delta(\pi) \rangle^2 &= \langle w^*, \Delta(\pi^*) \rangle^2 \\
&= \langle w^* , \Sigma w^* \rangle^2.
\end{align*}
\end{lem}
The proof of Lemma~\ref{lem:minimizing_u} is in Appendix~\ref{sec:proof:lem:minimizing_u}. 

By Proposition~\ref{prop:risk_wstar_pi}-(ii), we have  $\pi^* \in \Pi_{\epsilon,\xi,r} 
\subset \Pi_{\epsilon,\xi}$. It thus follows from Lemma~\ref{lem:minimizing_u} that
\begin{align}
   \inf_{\pi \in \Pi_{\epsilon,\xi,r}} u^2(w^*;\pi) & = 
   \inf_{\pi \in \Pi_{\epsilon,\xi}} \frac{\langle w^*, \Delta(\pi) \rangle^2}{\langle w^*, \Sigma w^* \rangle} = \langle w^* , \Sigma w^* \rangle \nonumber \\
   & = 2 N \sinh^2 \left( \frac{n \epsilon^2 N^{1-2/p} }{2} \right).
   \label{eq:inf_u}
\end{align}
We get
\begin{align}
    R^*(V_{\epsilon,\xi,r}) & = \inf_\psi \sup_{Q \in V_{\epsilon,\xi,r}} R(\psi,Q) \nonumber \\
    & \leq \sup_{Q \in V_{\epsilon,\xi,r}} R(\psi^*,Q) \nonumber \\
    & \leq \sup_{\pi \in \Pi_{\epsilon,\xi,r} } \rho(\psi^*;\pi) \label{eq:upper_proof_1},
\end{align}
the last transition 
because every $Q \in V_{\epsilon,\xi,r}$ corresponds to a prior in $\Pi_{\epsilon,\xi,r}$ given by the mass probability at $Q$. By \eqref{eq:upper_proof_1},
\begin{align*}
R^*(V_{\epsilon,\xi,r}) +o(1) & \leq \sup_{\pi \in \Pi_{\epsilon,\xi,r}} \rho(\psi^*;\pi) \\
& = \sup_{\pi \in \Pi_{\epsilon,\xi,r}}  2 \Phi(-u(w^*;\pi)/2) \\
& =  2 \Phi\left(- \inf_{\pi \in \Pi_{\epsilon,\xi,r}}  u(w^*;\pi)/2\right) 
\end{align*}
and \eqref{eq:prf:upper_bound:to_show} follows from \eqref{eq:inf_u}.
\qed

\subsection{Proof of Theorem~\ref{thm:lower_bound}
\label{sec:proof:lower_bound}}

Let $\xi = \log(N)/(n \sqrt{\log(\log(N))})$ and $r^2 = \log(N) / (n \sqrt{N})$. If $\log(N) / (n N^{1-1/p} \epsilon) = O(1)$, then $nN^{1-1/p} \epsilon \sqrt{\log(\log(N))} / 
\log(N) > 1$ for all $N$ large enough. For such $N$s, we have $N^{1-1/p} \epsilon > \xi$, thus $V_{\epsilon,\xi,r}=\emptyset$ by \eqref{eq:non_empty_constraint}, and thus $R^*(V_{\epsilon,\xi,r})=0$ by convention. In this case, Corollary~\ref{cor:trivial_separation} implies $R^*(V_\epsilon)=o(1)$. We henceforth assume that $\epsilon n N^{1-1/p} \leq \log(N)$ and thus $ \epsilon n N^{-1/p} = o(\log(N)/N)=o(1)$ and $\epsilon N^{1-1/p} = O(\log(N)/n)=o(1)$. 

From Lemma~\ref{lem:R_geq_rho}, we get
\begin{align}
    \inf_{\psi} \rho(\psi, \pi^*) & \leq  \rho^*\left(\Pi_{\epsilon,\xi,r} \right) = R^*(V_{\epsilon,\xi,r}) + o(1) \nonumber \\ 
    & \leq R^*(V_{\epsilon}) + o(1).
\label{eq:lower_bound_inequality_chain}
\end{align}
By Neyman-Pearson's theory, the test attaining the infimum in \eqref{eq:lower_bound_inequality_chain} is based on the likelihood ratio for testing the null hypothesis $H_0 = H(U)$ of \eqref{eq:hyp_Q} against a simple alternative of the form 
\begin{align*}
H_1\,& :\, O_i \sim \Pois(n Q_i),\qquad Q_i \sim \pi^*_i. 
\end{align*}
The log-likelihood ratio statistic is given by
\begin{align*}
\ell(O_1,\ldots,O_N;\pi^*) & := \sum_{i=1}^N \log(L_i(O_i;\pi_i^*)), 
\end{align*}
where 
\begin{align*}
L_i(x;\pi_i^*) & = \frac{\Prp{O_i=x \mid Q \sim \pi^*} }{\Prp{O_i = x \mid H_0} }.
\end{align*}
Note that,
\begin{align*}
L_i(x;\pi_i^*) & = \frac{1}{2}e^{n(1/N - \epsilon N^{-1/p})} \left( 1 + \epsilon N^{-1/p} \right)^{x} \\
& \qquad + \frac{1}{2}e^{n(\epsilon N^{-1/p} - 1/N)} \left( 1 - \epsilon N^{-1/p} \right)^{x} \\
& = 1 + g_{x,\lambda_0}(n \epsilon N^{-1/p}). 
\end{align*}
Using the identity
\[
\sum_{i=1}^N f(o_i) = \sum_{m=0}^\infty x_m f(m),\qquad x_m := \#\{i: o_i = m\},
\]
valid for an arbitrary function $f$, we conclude that the likelihood ratio test is of the form \eqref{eq:linear_test_stat} with weights
\begin{align}
\label{eq:LR_test_weights}
q_m= \log\left( 1 + g_{m,\lambda_0}(n \epsilon N^{-1/p}) \right).
\end{align}
Observe that if $\epsilon N^{1-1/p} < 1$, then $1 + g_{m,\lambda_0}\left( n \epsilon N^{-1/p} \right)>0$ ensuring that $q_m$ is well-defined.

We now show that the condition of Proposition~\ref{prop:asymp_normality} holds with $w_m = q_m$. Set $t_N := n \epsilon N^{-1/p}$. We assumed $t_N = o(1)$ and $\epsilon N^{1-1/p} = o(1)$, thus
 $t_N=o(\sqrt{\lambda_0})$. Provided $N$ is sufficiently large so $|t_N| \leq a < 1$, we may apply Lemmas~\ref{lem:gm_sums} and \ref{lem:gm_log_gm}. This gives 
\begin{align*}
    \sum_{m=0}^\infty \left|q_{m}\right|^4(t_N) \Po_{\lambda_0}(m) & = \left( \frac{t_N^2}{\lambda_0} \right)^4 \left(\frac{15}{4} + \frac{9}{\lambda_0} + \frac{1}{2 \lambda_0^2} + o(1) \right),
\end{align*}
\begin{align*}
    \left(\sum_{m=0}^\infty q_{m}(t_N) \Po_{\lambda_0}(m)\right)^2 & = o\left( \frac{t_N^2}{\lambda_0} \right)^2,
\end{align*}
and
\begin{align*}
    \sum_{m=0}^\infty q_{m}^2(t_N) \Po_{\lambda_0}(m) & = \frac{1}{2}\left( \frac{t_N^2}{\lambda_0} \right)^2 (1 + o(1)).
    \end{align*}
We use these bounds and $\sinh(x) = x(1+ o(1))$ for $x \to 0$ to obtain:
\begin{align}
    & \frac{\sum_{m=0}^\infty \left|\tilde{w}_m\right|^4 \Po_{\lambda_0}(m) }{N\left( \sum_{m=0}^\infty |\tilde{w}_m|^2 \Po_{\lambda_0}(m) - \left(\sum_{m=0}^\infty \tilde{w}_m \Po_{\lambda_0}(m) \right)^2 \right)^2} \nonumber \\
    & = \frac{1}{N} J\left( \frac{t_N^2}{\lambda_0} \right),
    \label{eq:log_m_w_cond}
\end{align}
where 
\begin{align*}
J(u) & := \frac{u^4(\frac{15}{4} + \frac{9}{\lambda_0} + \frac{1}{2 \lambda_0^2} +  o(1))}{\left( \frac{1}{2}u^2(1+o(1)) - o(u^2) \right)^2} \\
& = 4\left(\frac{15}{4} + \frac{9}{\lambda_0} + \frac{1}{2 \lambda_0^2} + o(1) \right) = O(N^2/n^2).
\end{align*}
By the assumption $N = o(n^2)$, \eqref{eq:log_m_w_cond} is $o(1)$ and \eqref{eq:w_cond2} follows. By Proposition~\ref{prop:asymp_normality},
\begin{align}
    \inf_{\psi} \rho(\psi,\pi^*) = 
\rho(\psi_{\tilde{w},\alpha^*},\pi^*) = 2 \Phi(-\tilde{u}_\epsilon/2) + o(1). 
\label{eq:prf:upper_LR_risk}
\end{align}
Combining \eqref{eq:prf:upper_LR_risk} with \eqref{eq:lower_bound_inequality_chain} yields \eqref{eq:prf:lower_bound}. \qed

\subsection{Proof of Theorem~\ref{thm:matching}
\label{sec:proof:matching}}
Set $u_\epsilon := u_{\epsilon,N,n,p}$ and $\xi$ and $r$ as in \eqref{eq:xi_r_explicit}. Consider here the limit $N \to \infty$ with $N=o(n^2)$. By Corollary~\ref{cor:trivial_separation}, the assumption $\lim R^*(V_\epsilon) > 0$ implies that $V_{\epsilon,\xi,r}$ is non-empty. By \eqref{eq:non_empty_constraint}, 
$\epsilon N^{1-1/p} \leq \log(N)/(n\sqrt{\log(\log(N))})$ for all $N$ large enough thus $\epsilon N^{1-1/p} = o(1)$. Likewise, $\epsilon N^{1/2-1/p} \leq r$, thus $n \epsilon N^{1-2/p} = O(\log(N)/\sqrt{N}) = o(1)$. 

Recall that $w_m^* = g_{m,\lambda_0}(n \epsilon N^{-1/p})$, $\tilde{w}_m = \log(1+ w_m^*)$, $\Delta_m(\pi^*) = w_m^* \Po_{\lambda_0}(m)$, 
and $\Sigma = \diag(\mu^0) - \mu^0 \mu^{0\top}/N$. 

By the inequality $\log(1+x) \leq x$ for $x>-1$, we get
\begin{align}
\langle \tilde{w},\Delta(\pi^*) \rangle & = N \sum_{m=0}^\infty \log (1 + w_m^*) w_m^* \Po_{\lambda_0}(m) \nonumber \\
& \leq N \sum_{m=0}^\infty (w_m^*)^2 \Po_{\lambda_0}(m) \nonumber \\
& = 2 N \sinh^2 \left( \frac{n \epsilon^2 N^{1-2/p}}{2} \right).
\label{eq:w_tilde_delta}
\end{align}
Additionally, 
\begin{align}
\langle \tilde{w}, \Sigma \tilde{w} \rangle & = \tilde{w}^\top \diag(\mu^0) \tilde{w} - (\tilde{w}^\top \mu^0)^2/N.
\label{eq:w_tilde_sigma_w_tilde}
\end{align}
Since $o(1) = n \epsilon^2 N^{1-2/p}$, it follows from Lemma~\ref{lem:gm_log_gm}-(iii) that
\begin{align}
    \label{eq:w_tilde_mu}
\left|\tilde{w}^\top \mu^0\right| = \left|\sum_{m=0}^\infty \tilde{w}_m \Po_{\lambda_0}(m) \right| \leq 2 \sinh^2 \left( \frac{n\epsilon^2 N^{1-2/p}}{2} \right) = o(1). 
\end{align}
This means that $(\tilde{w}^\top \mu^0)^2/N = o(1/N)$. By Lemma~\ref{lem:gm_log_gm}-(ii), 
\begin{align}
& \tilde{w}^\top \diag(\mu^0)\tilde{w} = N \sum_{m=0}^\infty \log^2\left(1+g_{m,\lambda_0}(n\epsilon N^{-1/p})\right) \Po_{\lambda_0}(m) \nonumber \\
& \quad \geq N \left( \sinh^2\left( \frac{n \epsilon^2 N^{1-2/p}}{2} \right) - 8\sinh^4\left( n \epsilon^2 N^{1-2/p} \right) \right) \nonumber \\
& \quad = N \sinh^2\left( \frac{n \epsilon^2 N^{1-2/p}}{2} \right)\left( 1 + o(1) \right).
\label{eq:w_tilde_w_tilde}
\end{align}
It follows from \eqref{eq:w_tilde_delta}, \eqref{eq:w_tilde_mu}, \eqref{eq:w_tilde_sigma_w_tilde}, and \eqref{eq:w_tilde_w_tilde} that 
\begin{align}
\tilde{u}_\epsilon^2 & = 
\frac{\langle \tilde{w},\Delta(\pi^*) \rangle^2}{\langle \tilde{w}, \Sigma \tilde{w} \rangle} \leq 2 N \sinh^2\left(\frac{n\epsilon^2 N^{1-2/p}}{2} \right) \left(1 + o(1) \right).
\label{eq:u_tilde_upper}
\end{align}
As $n \epsilon^2 N^{1-2/p} = o(1)$, we have 
\begin{align*}
2N \sinh^2\left(\frac{n\epsilon^2 N^{1-2/p}}{2} \right) & = 2N\left( n\epsilon^2 N^{1-2/p}/2 \right)^2(1+o(1)) \\
& = {u_\epsilon^2}(1+o(1)). 
\end{align*}
Therefore, $\tilde{u}_\epsilon^2 \leq u_\epsilon^2(1+o(1))$.
By Theorems~\ref{thm:upper_bound} and \ref{thm:lower_bound}, and uniform continuity of the normal CDF,
\begin{align*}
 R^*(V_\epsilon) + o(1) & \geq 2\Phi(-\tilde{u}_\epsilon/2) + o(1) \\
 & \geq 2\Phi\left(- u_{\epsilon}(1+o(1))/2
 \right)+ o(1) \\
 & =  2\Phi\left(- u_{\epsilon}/2
 \right) + o(1) \geq R^*(V_\epsilon) + o(1),
\end{align*}
thus 
\[
R^*(V_\epsilon) + o(1) = 2\Phi\left(- u_{\epsilon}/2
 \right). 
\]
From here, the assumption $\lim R^*(V_\epsilon) \in (0,1]$ implies $\lim u_\epsilon \in [0,\infty)$. This completes the proof. \qed

\subsection*{Proof of Proposition~\ref{prop:chisq_n_collision}
\label{sec:chisq_n_collision:proof}}
With $w_m = (m-\lambda_0)^2/\lambda_0$, by formulas of the centralized Poisson moments we get $(\mu^0)\top w = N$ and $w^\top \diag(\mu^0)w=N(3 + 1/\lambda_0)$, hence
\[
w^\top \Sigma w = N(3 + \frac{1}{\lambda_0} - 1) = N(2 + \frac{1}{\lambda_0}).
\]
In addition, 
\begin{align*}
w^\top \Delta(\pi^*) & = N \sum_{m=0}^\infty \Po_{\lambda_0}(m) \frac{(m-\lambda_0)^2}{\lambda_0} g_{m,\lambda_0}(n \epsilon N^{-1/p}) \\
& = N \frac{(n \epsilon N^{-1/p})^2}{\lambda_0} = \frac{n^2 N^{2-2/p} \epsilon^2}{\lambda_0}.
\end{align*}
It follows that 
\begin{align*}
 \frac{\langle w , \Delta(\pi^*) \rangle^2}{\langle w, \Sigma w \rangle }&= \frac{n^2 N^{3-4/p}\epsilon^4/\lambda_0^2}{N(2 + 1/\lambda_0)} \\
&= \frac{n^2 N^{2-4/p}\epsilon^4}{\lambda_0^2(2 + 1/\lambda_0)} \\
&=  u_{\epsilon,n,N,p}^2\frac{\lambda_0}{2\lambda_0 + 1}.
\end{align*}
The proof is concluded by showing that the conditions of Proposition~\ref{prop:asymp_normality} are satisfied. Condition \eqref{eq:w_cond1} holds since the second moment of the Poisson distribution exists. For Condition \eqref{eq:w_cond2}, consider the following evaluations:
\begin{align*}
    \sum_{m=0}^\infty w_m \Po_{\lambda_0}(m) & = \sum_{m=0}^\infty \frac{(m-\lambda_0)^2}{\lambda_0} \Po_{\lambda_0}(m) = 1,
\end{align*}
\begin{align*}
    \sum_{m=0}^\infty \left|w_m\right|^2 \Po_{\lambda_0}(m) & = \sum_{m=0}^\infty \frac{(m-\lambda_0)^4}{\lambda_0^2} \Po_{\lambda_0}(m) \\
    & = \frac{\lambda_0 + 3\lambda_0^2}{\lambda_0^2} = \frac{1 + 3\lambda_0}{\lambda_0},
\end{align*}
and
\begin{align*}
    \sum_{m=0}^\infty \left|w_m\right|^4 \Po_{\lambda_0}(m) & = \sum_{m=0}^\infty \frac{(m-\lambda_0)^8}{\lambda_0^4} \Po_{\lambda_0}(m) \\
    & = \frac{1}{\lambda_0^4} \left( 105 \lambda_0^4 + 490 \lambda_0^3 + 119 \lambda_0^2 + \lambda_0 \right).
\end{align*} 
We have
\begin{align*}
    & \frac{\sum_{m=0}^\infty \left|w_m\right|^4 \Po_{\lambda_0}(m)}{N \left( \sum_{m=0}^\infty \left|w_m\right|^2 \Po_{\lambda_0}(m)  - \left(\sum_{m=0}^\infty w_m \Po_{\lambda_0}(m) \right)^2 \right)^2 }\\
    & = \frac{105 \lambda_0^4 + 490 \lambda_0^3 + 119 \lambda_0^2 + \lambda_0}{N \lambda_0^2(1 + 2\lambda_0)^2}
\end{align*}
The last expression is $O(1/N)$ when $\lambda_0$ is bounded from below and $O(1/n)$ when $\lambda_0 = o(1)$. In any case, condition \eqref{eq:w_cond2} is satisfied. This proves \eqref{eq:asymp_risk_chisq}.
\\

For the collision statistic, we have $w_m = \binom{m}{2} = \frac{m(m-1)}{2}$. By standard Poisson moments, $(\mu^0)^\top w = N \cdot \frac{\lambda_0^2}{2}$ and $w^\top \diag(\mu^0) w = N \left(\lambda_0^3 + \frac{\lambda_0^2}{2}\right)$. Therefore,
\[
w^\top \Sigma w = N \left(\lambda_0^3 + \frac{\lambda_0^2}{2}\right) = N \frac{\lambda_0^2(2\lambda_0+1)}{2}. 
\]
Additionally, 
\[
w^\top \Delta(\pi^*) = N \ex{\binom{M}{2} g_{M,\lambda_0}(x)}.
\]
Let $M \sim \Pois(\lambda_0)$. Using the Poisson probability generating function identity,
\[
\ex{\binom{M}{2} t^M} = \frac{1}{2} t^2 \lambda_0^2 e^{\lambda_0(t-1)},
\]
with $t_\pm = 1 \pm x/\lambda_0$, we get
\begin{align*}
\ex{\binom{M}{2} g_{M,\lambda_0}(x)}
&= \frac12\Big(
e^{-x}\ex{\binom{M}{2} t_+^M}
+ e^{x}\ex{\binom{M}{2} t_-^M} \\
& \qquad - 2\ex{\binom{M}{2}}
\Big) \\
&= \frac12\Big(
\frac{\lambda_0^2}{2}(t_+^2+t_-^2) - \lambda_0^2 \Big)
= \frac{x^2}{2}.
\end{align*}
Hence, with $x = n\epsilon N^{-1/p}$,
\begin{align*}
w^\top \Delta(\pi^*)
& = N \frac{x^2}{2} = \frac{N}{2}(n\epsilon N^{-1/p})^2 \\
& = \frac{n^2\epsilon^2}{2} N^{1-2/p}.
\end{align*}
Finally,
\begin{align*}
\frac{\langle w , \Delta(\pi^*) \rangle^2}{\langle w, \Sigma w \rangle}
& = \frac{\left(\frac{n^2\epsilon^2}{2} N^{1-2/p}\right)^2}
{N\frac{\lambda_0^2(2\lambda_0+1)}{2}} \\
& = \frac{n^4\epsilon^4}{2\lambda_0^2(2\lambda_0+1)}\,N^{1-4/p} \\
& = u_{\epsilon,n,N,p}^2\frac{1}{2\lambda_0+1}.
\end{align*}

Verification that the conditions of Proposition~\ref{prop:asymp_normality}
hold for these weights proceeds similarly to the case $w_m = (m-\lambda_0)^2/\lambda_0$ for the chi-squared weights and is omitted.
This proves \eqref{eq:asymp_risk_col}. \qed

\appendices



\section{Proofs of Technical Lemmas
\label{app:proof_of_lemmas}
}
This Appendix contains the proofs of several technical lemmas used in the proof of the main results in Section~\ref{sec:proofs} above. 

\subsection{Proof of Lemma~\ref{lem:asym_normality1}
\label{sec:proof:lem:asym_normality1}
}
Consider
\[
T_w = \sum_{m=0}^\infty w_m X_m =  \sum_{i=1}^N w_{O_i},\quad O_i \sim \Pois(nQ_i). 
\]
We will show that the sequence of random variables $\{w_{O_i}\}_{i=1}^N$ satisfies the condition for a Berry-Esseen type central limit theorem (c.f. \cite[Appendix A]{shorack2009empirical}):
\begin{align}
    E_n & := 
    \frac{ \sum_{i=1}^N \ex{\left| w_{O_i}\right|^4} }{ \left(\sum_{i=1}^N \Var{w_{O_i}} \right)^{2}} \to 0,
    \label{eq:Berry-Esseen-cond}
\end{align}
as $N\to \infty$. Notice that $\{w_m\}$, $\lambda$, and $Q$ generally depend on $N$ and $n$. For $v >0$, denote $\Upsilon_v \sim \Pois(v)$. Under the null $Q=U$, 
\begin{align*}
\Var{w_{O_i}} & = \Var{w_{\Upsilon_\lambda}} \\
& = \sum_{m=0}^\infty w_m^2 \Po_\lambda(m) - \left( \sum_{m=0}^\infty w_m \Po_\lambda(m) \right)^2 \\
& = \sum_{m=0}^\infty \left(\sum_{k=0}^\infty w_k \Po_{\lambda}(k) -w_m \right)^2 \Po_\lambda(m).
\end{align*}
Because $\{w_m\}$ is not a constant sequence, there exists $m$ such that
$\sum_{k=0}^\infty w_k \Po_{\lambda}(k) \neq w_m$ and thus $\Var{w_{\Upsilon_\lambda}} >0$. For $a \geq 2$, denote
\[
\gamma_{a}(v) := \ex{\left|w_{\Upsilon_{v}}\right|^a} = \sum_{m=0}^\infty \left| w_m \right|^a \Po_{v}(m). 
\]
(condition \eqref{eq:w_cond1} ensures that the sum is well-defined for any $v>0$.) Under the null, $\ex{\left|w_{O_i} \right|^4} = \gamma_{4}(\lambda)$. Therefore, under the null,
\begin{align}
\label{eq:E_n0}
E_n & =  \frac{N \gamma_{4}(\lambda)}{\left( N \Var{w_{\Upsilon_{\lambda}}} \right)^{2}} \\
& = \frac{\sum_{m=0}^\infty \left|w_m\right|^4 \Po_{\lambda}(m) }{N\left( \sum_{m=0}^\infty |w_m|^2 \Po_{\lambda}(m) \right.} \nonumber \\
& \qquad \left. - \left(\sum_{m=0}^\infty w_m \Po_{\lambda}(m) \right)^2 \right)^2.
\end{align}
From here, \eqref{eq:w_cond2} implies \eqref{eq:Berry-Esseen-cond} which implies
\eqref{eq:lem1:statement}. 

Henceforth, we assume $\Var{w_{\Upsilon_{\lambda}}} = 1$. This assumption is without loss of generality because for $\tilde{w}_m := w_m/\sqrt{\Var{w_{\Upsilon_{\lambda}}}}$, we have
\begin{align*}
\tilde{E}_n &:= 
    \frac{ \sum_{i=1}^N \ex{\left| \tilde{w}_{O_i}\right|^4} }{ \left(\sum_{i=1}^N \Var{\tilde{w}_{O_i}} \right)^{2}} \\
&= \frac{ \sum_{i=1}^N \ex{\left| w_{O_i}\right|^4} }{ \left(\sum_{i=1}^N \Var{w_{O_i}} \right)^{2}} = E_n.
\end{align*}
In particular, we have $1\leq \gamma_{4}(\lambda)$. 

By \eqref{eq:w_cond1}, 
\begin{align}
\ex{\left|w_{\Upsilon_{\lambda}}\right|^a} & = \sum_{m=0}^\infty \left|w_m\right|^a \Po_{\lambda}(m) \nonumber  \leq \sum_{m=0}^\infty C_0e^{C_1 \cdot m\cdot a} \Po_{\lambda}(m) \\
& = C_0 e^{\lambda(e^{aC_1}-1)} = C_0 e^{2\lambda C_a},
\label{eq:moment_bound}
\end{align}
for some constant $C_a$ that only depends on $a$. 

For two distributions over $\{0,1,2,\ldots\}$ $P_1$ and $P_2$, let $\mathrm{TV}(P_1,P_2)= \frac{1}{2}\sum_{m=0}^\infty \left| P_1(m)-P_2(m)\right|$. Suppose that 
\[
\sum_{m=0}^\infty \left|v_m \right|^{2a} P_j(m) < \infty,\qquad j=1,2,
\]
for some sequence $\{v_m\}$. Then by Cauchy-Schwarz inequality,
\begin{align*}
    & \left( \sum_{m=0}^\infty \left|v_m\right|^a \left|P_1(m)-P_2(m)\right| \right)^2 \\
    & \leq \sum_{m=0}^\infty|v_m|^{2a} \left|P_1(m)-P_2(m)\right| \sum_{m=0}^\infty \left|P_1(m)-P_2(m)\right| \\
    & \leq \left[\sum_{m=0}^\infty |v_m|^{2a} P_1(m) + \sum_{m=0}^\infty |v_m|^{2a} P_2(m) \right] \cdot 2 \mathrm{TV}(P_1,P_2). 
\end{align*}
Therefore, using that $\mathrm{TV}(\Upsilon_{\lambda_1}, \Upsilon_{\lambda_2}) \leq |\lambda_1 - \lambda_2|$ (c.f. 
\cite[Prop. 5]{kelbert2023survey}), 
\begin{align*}
    \left(\gamma_a(\lambda_1) - \gamma_a(\lambda_2)\right)^2 & = \left(\sum_{m=0}^\infty 
    |w_m|^a \left(\Po_{\lambda_1}(m) - \Po_{\lambda_2}(m) \right) \right)^2 \\
    & \leq 
    \left(\sum_{m=0}^\infty 
    |w_m|^a \left|\Po_{\lambda_1}(m) - \Po_{\lambda_2}(m) \right|\right)^2
    \\
    & \leq 
    2\left[\sum_{m=0}^\infty |w_m|^{2a} \Po_{\lambda_1}(m) \right. \nonumber \\
    & \qquad \left. + \sum_{m=0}^\infty |w_m|^{2a} \Po_{\lambda_2}(m) \right] \left|\lambda_1-\lambda_2\right|.
\end{align*}
From \eqref{eq:moment_bound} and since $\max_i |Q_i| \leq b_{N}$ and $1/N \leq b_{N}$, 
we get the following by taking $\lambda_1 = n/N$ and $\lambda_2 = n Q_i$.
\begin{align*}
\left|\gamma_a(nQ_i) - \gamma_a(\lambda_0)\right|^2 
    & \leq  4C_0 e^{n b_{N} C_{2a}} \left|nQ_i-n/N \right|. 
\end{align*}
By Jensen's inequality and concavity of $x \to x^{1/4}$, we get
\begin{align*}
& \sum_{i=1}^N \left|\gamma_a(nQ_i) - \gamma_a(\lambda_0) \right| \leq \sqrt{4C_0} e^{ n b_{N} C_{2a}} 
\sum_{i=1}^N  \left( |nQ_i-n/N|^2 \right)^{1/4} \\
& \quad \leq 
 \sqrt{4C_0} e^{  n b_{N} C_{2a}} \left( n^2 N^3 \sum_{i=1}^N |Q_i-1/N|^2 \right)^{1/4} \\
 & \quad = \sqrt{4C_0} e^{ n b_{N} C_{2a}} \left( n^2 N^3 \left\| Q - U \right\|_2^2 \right)^{1/4}.
\end{align*}
Consequently, by \eqref{eq:Q_cond},  
\begin{align*}
\sum_{i=1}^N \left|\gamma_a(nQ_i) - \gamma_a(\lambda_0) \right| &= \sqrt{4C_0} e^{n b_{N} C_{2a}}\left( n^2 N^3 a_{N} \right)^{1/4} \\
&= o(N), 
\end{align*}
where the last transition is because $n^2 N^{-1} a_{N} e^{C_{2 a} n b_{N}} =o(1)$ by assumption. It follows that 
\[
\left|\sum_{i=1}^N \gamma_a(n Q_i) - \sum_{i=1}^N \gamma_a(\lambda_0)\right| = o(N),
\]
hence,  
\begin{align}
& \frac{
\sum_{i=1}^N \gamma_a(n Q_i)}
{\sum_{i=1}^N \gamma_a(\lambda_0)} = \frac{
\sum_{i=1}^N \gamma_a(n Q_i)}
{N \gamma_a(\lambda_0)} \nonumber \\
& \qquad = 1 + \frac{o(N)}{N \gamma_a(\lambda_0)} = 1 + o(1/\gamma_a(\lambda_0)).
\label{eq:H1_H0_gamma_equivalnce}
\end{align}
By similar arguments,
\begin{align}
\left|\ex{w_{O_i}} -  \ex{w_{\Upsilon_{\lambda_0} }} \right| \leq \ex{w_{\Upsilon_{\lambda_0} }} + o(1).
\label{eq:mean_equivalence}
\end{align}
From \eqref{eq:H1_H0_gamma_equivalnce}  and the assumption $\Var{w_{\Upsilon_{\lambda_0}}}=1$ (thus $1 \leq \gamma_4(\lambda_0)$), we obtain 
\[
\frac{
\sum_{i=1}^N \gamma_4(n Q_i)}
{\sum_{i=1}^N \gamma_4(\lambda_0)} = 1 + o(1),
\]
and, with the help of \eqref{eq:mean_equivalence},
\begin{align}
\frac{\sum_{i=1}^N \Var{w_{O_i}}}{\sum_{i=1}^N \Var{w_{\Upsilon_{\lambda_0}}}} &= \frac{\sum_{i=1}^N \Var{w_{O_i}}}{N} \nonumber \\
&= 1 + o(1).
\label{eq:variance_equivalence}
\end{align}
It follows that
\begin{align*}
E_n & = \frac{ \sum_{i=1}^N \gamma_{4}(nQ_i)}{ ( \sum_{i=1}^N \Var{w_{O_i}} )^{2}} \\
&=  \frac{(1+o(1))N \gamma_{4}(\lambda_0)}{(N (1+o(1)))^{2}} = o(1);
\end{align*}
the last transition is due to \eqref{eq:E_n0} and $\Var{w_{\Upsilon_{\lambda_0}}} = 1$. The Berry-Esseen condition \eqref{eq:Berry-Esseen-cond} applies, leading to
\begin{align*}
\lim_{n \to \infty} \sup_{x\in \reals} \left[\Prp{ \frac{T_w - \ex{T_w}}{\sqrt{\Var{T_w}}} \leq x \mid H(Q)} - \Phi(x) \right] =0.
\end{align*}
Finally, \eqref{eq:variance_equivalence}  says
\begin{align*}
    \Var{T_w} & = 
    \sum_{i=1}^N \Var{w_{O_i}} \\
    &= (1+o(1)) \sum_{i=1}^N \Var{w_{\Upsilon_{\lambda_0}}} \\
    & = (1+o(1))\Var{T_w|H(U)} \\
    &= (1+o(1))\langle w, \Sigma w \rangle,
\end{align*}
hence 
\begin{align*}
\lim_{n \to \infty} \sup_{x\in \reals} \left[\Prp{ \frac{T_w - \ex{T_w}}{\sqrt{\langle w, \Sigma w \rangle}} \leq x \mid H(Q)} - \Phi(x) \right] =0.
\end{align*} 
\qed

\subsection{Proof of Lemma~\ref{lem:minimizing_u}
\label{sec:proof:lem:minimizing_u}}

Denote
\[
G(\pi) := \langle w^*, \Delta(\pi) \rangle^2. 
\]
We have
\begin{align}
    \label{eq:wstar_delta}
& \langle w^*, \Delta(\pi) \rangle  \\
& = \sum_{m=0}^\infty w^*_m \Po_{\lambda_0}(m)  \sum_{i=1}^N \int_{\reals} h_{m,\lambda_0}(n(t-1/N))\pi_i(dt). 
\end{align}
To simplify notation, we implicitly assume throughout this proof that the $i$-th coordinate of a prior $\pi$ is shifted by $1/N$. For example, instead of \eqref{eq:wstar_delta} we write 
\begin{align*}
& \langle w^*, \Delta(\pi) \rangle \\
& \quad = \sum_{m=0}^\infty w^*_m \Po_{\lambda_0}(m)  \sum_{i=1}^N \int_{\reals} h_{m,\lambda_0}(nt)\pi_i(dt). 
\end{align*}
We also denote $\Pi = \Pi_{\epsilon,\xi}$ to reduce clutter.

We first reduce attention to the subset of measures $\bar{\Pi} \subset \Pi$ in which each coordinate is symmetric around $0$ (in the original notation, symmetric around $1/N$). For a one-dimensional measure $\pi_1$ on $\reals$, denote $\pi_1^{\#}(dt) = \pi_1(-dt)$ and 
\[
\bar{\pi}_1(dt) = \frac{\pi_1(dt) + \pi_1^{\#}(dt)}{2},
\]
and by $\bar{\pi}$ and $\pi^{\#}$ the obvious extension of these operations to product priors. Notice that if $\pi \in \Pi$, then $\bar{\pi} \in \Pi$, i.e. $\bar{\Pi} \subset \Pi$. Since $\langle w^*, \Delta(\pi) \rangle$ is linear in $\pi$, the function $G(\pi)$ is convex and hence
\begin{align}
    \label{eq:proof:pi_star:symmetric}
G(\bar{\pi}) \leq \frac{1}{2} G(\pi) + \frac{1}{2} G(\pi^{\#}).
\end{align}
Because symmetry dictates
\[
\inf_{\pi \in \Pi} G(\pi) = \inf_{\pi \in \Pi} G(\pi^{\#}), 
\]
it follows from \eqref{eq:proof:pi_star:symmetric} that
\[
\inf_{\pi \in \Pi} G(\bar{\pi}) \leq \inf_{\pi \in \Pi} G(\pi). 
\]
Thus, it is enough to consider the minimum of $G(\pi)$ over $\bar{\Pi}$. We identify each $\bar{\pi} \in \bar{\Pi}$ with a one-sided prior $\pi$ 
with support in $[0,\xi]$ such that 
\[
2 \sum_{i=1}^N \int_{\reals_+} |t|^p \pi_i(dt) \geq  \epsilon^p. 
\]
We denote the set of priors defined this way by $\Pi^+$. Because $2g_{m,\lambda}(x):= h_{m,\lambda}(x) + h_{m,\lambda}(-x)$, for $\bar{\pi} \in \Pi^+$ we have
\begin{align*}
    \Delta_m(\bar{\pi}) = \Po_{\lambda}(m) \sum_{i=1}^N \int_{\reals_+} g_{m,\lambda}(nt) \pi_i(dt),\qquad m=0,1,...
\end{align*}
Consider the formal sum 
\begin{align*}
f(t) &:= g(t;n,N,p,\epsilon) \\
&:= \sum_{m=0}^\infty \Po_{\lambda_0}(m)  g_{m,\lambda_0}\left(nN^{-1/p}\epsilon \right) g_{m,\lambda_0}(nt).
\end{align*}
Since all moments of the Poisson distribution exist, this sum defines a smooth function $f(t)$. Furthermore, 
\begin{align*}
f(t) &= \cosh(N^{1-1/p}\epsilon n t)-1 \\
&= 2\sinh^2\left(N^{1-1/p}\epsilon n t/2\right), 
\end{align*}
as can be verified by standard power series identities and Taylor expansion of the function $\cosh(x)$. 

For a one-dimensional prior $\pi_i$ over $\reals^+ = [0,\infty)$, we set 
\begin{align*}
    g(\pi_i) := \exsub{Q_i \sim \pi_i}{f(Q_i)} = \int_{\reals_+} f(t)  \pi_i(dt). 
\end{align*}
For $\pi \in \Pi^+$, 
\begin{align}
G(\pi) & = \left(\langle w^*, \Delta(\pi) \rangle \right)^2 \nonumber \\
& = \left(\sum_{m=0}^\infty w_m^* \Po_{\lambda_0}(m) \sum_{i=1}^N  \int_{\reals_+} g_{m,\lambda_0}(nt) \pi_i(dt) \right)^2 \nonumber \\
& = \left(\sum_{i=1}^N \int_{\reals_+} f(t) \pi_i(dt) \right)^2 \nonumber \\
& = \left(\sum_{i=1}^N g(\pi_i) \right)^2 
\label{eq:proof:pi_i}
\end{align}
Consider the set of sequences
    \[
    L_{\epsilon} := \{ (a) \in \reals^N \,:\, \sum_{i=1}^N  a_i^p \geq \epsilon^p,\quad 0 \leq a_i\},
    \] 
and the set of one-dimensional positive priors 
\[
\Pi_{a}^+ = \{ \pi_1 \,:\, \int_{\reals_+} t^p \pi_1 (dt) \leq a^p,\, \pi_1([0, \xi])=1 \}.
\]
We have
\begin{align}
\inf_{\pi \in \Pi^+} G(\pi) & = \left( \inf_{\pi \in \Pi^+} \sum_{i=1}^N g(\pi_i) \right)^2 \nonumber \\
& = \inf_{(a) \in L_{\epsilon} } \left( \sum_{i=1}^N \inf_{\pi_i \in \Pi_{a_i}^+ } g(\pi_i)  \right)^2 \label{eq:proof:two_inf} 
 \\ 
 & =: \inf_{(a) \in L_{\epsilon} } \left( \sum_{i=1}^N g_{a_i} \right)^2,
 \nonumber
\end{align}
where we denoted $g_a := \inf_{\pi_1 \in \Pi^+_a} g(\pi_1)$. 

The following lemma addresses the inner minimization in \eqref{eq:proof:two_inf}. 
\begin{lem}
\label{lem:three_point}
Fix $a>0$. There exist $\eta_{a} \in [0,1]$ and $\mu_{a} \geq0$ such that, for $\pi^*_1 = (1-\eta_{a}) \delta_0 + \eta_{a} \delta_{\mu_a} \in \Pi^+_{a}$, 
\begin{align*}
g_a &= \inf_{\pi_1 \in \Pi^+_a} g(\pi_1) \\
&= g(\pi_1^*) = \eta_a f(\mu_a).
\end{align*}
\end{lem}
Namely, the inner minimum in \eqref{eq:proof:two_inf} is attained by non-negative two-point priors. Such a prior is later identified with a symmetric three-point prior (or a two-point prior when $\eta_t = 1$). The proof of Lemma~\ref{lem:three_point} is in Appendix~\ref{sec:proof:lem:three_points}. 

For the outer minimization in \eqref{eq:proof:two_inf}, we use: 
\begin{lem}
    \label{lem:equal_lambda}
Let $\phi$ be a real, convex, and permutation invariant function over $N$ variables. Then
\[
\phi(x_1,\ldots,x_N) \geq \phi\left(\bar{x},...,\bar{x} \right), \qquad \bar{x} := \frac{1}{N}\sum_{i=1}^N x_i.
\]
\end{lem}
The proof of Lemma~\ref{lem:equal_lambda} is in Appendix~\ref{sec:proof:lem:equal_lambda}. Applying Lemma~\ref{lem:equal_lambda} to the function 
\begin{align*}
\phi(x_1,\ldots,x_N) = \left(\sum_{i=1}^N x_i \right)^2,
\end{align*}
we see that the minimum in the outer minimization in \eqref{eq:proof:two_inf} is attained at sequences $(a)$ with $a_1=\ldots=a_N$. Since $g_{a_1}$ is increasing in $a_1 > 0$, for such a minimal-attaining sequence the constraint on $(a)$ is attained with equality, leading to  $\epsilon^p = N a_1^p$. Set $b = \epsilon N^{-1/p}$. Parametrize the set $\bar{\Pi}^3_{b}$ of one-dimensional symmetric three-point priors by $(\eta, \mu)$, so that for $\pi_1 \in \bar{\Pi}^3_{b}$ we have $g(\pi_1) = \eta f(\mu)$. The right-hand side of \eqref{eq:proof:two_inf} now leads to 
\begin{align}
\inf_{\pi \in \Pi^+} G(\pi) & = \left( \sum_{i=1}^N \inf_{\pi_i \in \Pi^3_b} g(\pi_i) \right)^2 \label{eq:proof:two_inf3} \\
& = \left( \sum_{i=1}^N \inf_{\pi_1 \in \Pi^3_b} g(\pi_1) \right)^2 \label{eq:proof:two_inf1} \\
& = N^2 \inf_{\eta,\mu} \eta^2 f^2(\mu), \label{eq:proof:two_inf4}
\end{align}
where the last minimization is over $\eta$ and $\mu$ such that 
\begin{align*}
(1-\eta)\delta_0 + \frac{\eta}{2}\delta_{\mu} + \frac{\eta}{2}\delta_{-\mu} \in \Pi^3_b,\\
b = \epsilon N^{-1/p}.
\end{align*}
Note that \eqref{eq:proof:two_inf3} follows from Lemma~\ref{lem:three_point}  and \eqref{eq:proof:two_inf1} follows from Lemma~\ref{lem:equal_lambda} as explained earlier. Equation \eqref{eq:proof:two_inf4} describes a one-dimensional problem that is equivalent to 
\begin{align}
\begin{split}
    \label{eq:optimization_problem}
    \text{minimize:~~} & \eta^2 f^2(\mu) \\
    \text{subject to:~~} & \eta N\mu^p \geq \epsilon^p,\quad 0 \leq \eta \leq 1 ,\quad 0\leq \mu.
    \end{split}
\end{align}
For $p \leq 2$, the solution is $\eta^* = 1$ and $\mu^* = \epsilon N^{-1/p}$, obtained by the method of Lagrange multipliers and convexity of $f(t)$. This completes the proof. 
\qed

\subsection{Proof of Lemma~\ref{lem:three_point}
\label{sec:proof:lem:three_points}
}
The proof uses similar ideas as the proof of \cite[Lemma 2]{ingster1994minimax} and \cite[Lemma 1]{suslina1999extremum}.  
For $\pi_1 \in \Pi_{a}^+$, define the measure $r(dt)=\pi_1(dt) - \pi_1(\{0\})\delta_0$ and the measure $q_1(dt)$ by $d q_1/d r = a^{-p} t^p$. We have
\[
q_1(\reals_+) = q_1([0,\xi]) = 1\quad \text{and} \quad \int t^{-p} q_1(dt) \leq a^{-p}. 
\]
We can write the minimization in question,
\begin{align}
g_{a} = \inf_{\pi_1 \in \Pi^+_{a}} g(\pi_1) = \inf_{\pi_1 \in \Pi^+_{a}} \int_{\reals_+} f(s)  \pi_1(ds),
    \label{eq:lem:three_points}
\end{align}
as 
\begin{align}
\label{eq:proof:two_inf_lemma}
    \inf_{q_1} \left\{ \int_{\reals_+} f(t) q_1(dt),\,:\, q_1 ([0,\xi])=1, \right. \nonumber \\
    \left. \int_{\reals_+} t^{-p} q_1(dt) \leq a^{-p}
    \right\}
\end{align}
By convexity of $f(s)=2\sinh^2(s)$ and Jensen's inequality, 
\begin{align*}
     \int_{\reals_+} f(s) q_1(ds) & \geq f\left( \int_{\reals_+} s q_1(ds) \right). 
\end{align*}
Likewise, convexity of $t \to t^{-p}$ implies 
\[
\left( \int_{\reals_+} t q_1(dt) \right)^{-p} \leq \int_{\reals_+} t^{-p} q_1(dt) \leq a^{-p}. 
\]
The last two inequalities prove that the minimum in \eqref{eq:proof:two_inf_lemma} is attained when $q_1= \delta_{a}$, hence the minimum in \eqref{eq:lem:three_points} is attained by $\pi_1 = (1-\eta)\delta_0 + \eta \delta_{a}$ for $\eta \geq 0$. 
\qed

\subsection{Proof of Lemma~\ref{lem:equal_lambda}
\label{sec:proof:lem:equal_lambda}
}

Let $\sigma$ indicate a uniform random variable over the permutation set of $\{1,...,N\}$. By permutation invariance,
\begin{align*}
    \phi(x_1,...,x_N) & = \phi( x_{\sigma(1)},...,x_{\sigma(N)}),
\end{align*}
for any $x_1,...,x_N$ in the domain of $\phi$, hence $\phi(x_1,...,x_N) = \ex{\phi\left( x_{\sigma(1)},...,x_{\sigma(N)} \right)}$. Therefore, by Jensen's inequality
\begin{align*}
    \phi(x_1,...,x_N) & = \ex{\phi\left( x_{\sigma(1)},...,x_{\sigma(N)} \right)} \\
    & \geq \phi \left( \ex{x_{\sigma(1)},...,x_{\sigma(N)}} \right)  \\
    & = \phi\left(\bar{x},\ldots,\bar{x} \right). 
\end{align*}
\qed

\section{Properties of \texorpdfstring{$h_{m,\lambda}$ and $g_{m,\lambda}$}{h and g}}\label{app:gm}

In this Appendix, we prove several properties of the functions 
\[
h_{m,\lambda}(x) = e^{-x}(1+\frac{x}{\lambda})^m - 1,
\]
and 
\[
g_{m,\lambda}(x) = \frac{h_{m,\lambda}(x) + h_{m,\lambda}(-x)}{2}
\]
that were used in some of the proofs in Section~\ref{sec:proofs} and Appendix~\ref{app:proof_of_lemmas}.

\begin{lem}
\label{lem:gm_properties}
Consider the function $g_{m,\lambda}(x)$, $\lambda>0$, and $m=0,1,...$. The following holds. 
\begin{enumerate}
\item[(i)] For any $x \in \reals$ and $\lambda'>0$, 
\begin{align*}
\sum_{m=0}^\infty g_{m,\lambda'}(x) \Po_{\lambda}(m)= 2 \sinh^2\left(\frac{x}{2}\left( \frac{\lambda}{\lambda'}-1\right) \right).
\end{align*}
\item[(ii)] For any $x \in \reals$ and $\lambda'>0$, 
\begin{align*}
\sum_{m=0}^\infty m g_{m,\lambda'}(x) \Po_{\lambda}(m) &= 2 \lambda \sinh^2\left( \frac{x}{2} \left(\frac{\lambda}{\lambda'}-1\right) \right) \\
&\quad + x\frac{\lambda}{\lambda'} \sinh \left( x \left(\frac{\lambda}{\lambda'}-1\right) \right).
\end{align*}
\item[(iii)] For any $t \in \reals$ and $m=0,1,...$, 
\begin{align}
\label{eq:gm_appx}
    g_{m,\lambda}(t \lambda) & = \frac{t^2}{2}\left( m(m-1) - 2 \lambda m + \lambda^2 \right)+o(t^2 (e^m + \lambda^2)).
    \end{align}
\item[(iv)] Fix $c \in (0,1)$. For $t \in [-c\lambda, c \lambda]$, we have $-1 < (1-c)^m/2 - 1 \leq g_{m,\lambda}(t)$. 
\end{enumerate}
\end{lem}

\subsection*{Proof of Lemma~\ref{lem:gm_properties}}
We have
\begin{align*}
    \sum_{m=0}^\infty e^{-x}\left(1+\frac{x}{\lambda'}\right)^m e^{-\lambda} \frac{\lambda^m}{m!} = e^{x\left(\frac{\lambda}{\lambda'}-1 \right)}.
\end{align*}
Thus, (i) follows from
\begin{align*}
\sum_{m=0}^\infty g_{m,\lambda'}(x) \Po_{\lambda}(m) & = \frac{1}{2}e^{x\left(\frac{\lambda}{\lambda'}-1\right)} + \frac{1}{2}e^{-x\left(\frac{\lambda}{\lambda'}-1\right)} - 1 \\
& = \cosh\left(x\left( \frac{\lambda}{\lambda'}-1\right) \right) - 1  \\
& = 2\sinh^2\left(\frac{x}{2}\left( \frac{\lambda}{\lambda'}-1\right) \right).
\end{align*}
Likewise, 
\begin{align*}
    \sum_{m=0}^\infty m e^{-x}\left(1+\frac{x}{\lambda'}\right)^m e^{-\lambda} \frac{\lambda^m}{m!} = e^{x\left(\frac{\lambda}{\lambda'}-1 \right)}\left(1 + \frac{x}{\lambda'} \right)\lambda,
\end{align*}
so 
\begin{align*}
\sum_{m=0}^\infty m g_{m,\lambda'}(x) \Po_{\lambda}(m) &= \frac{\lambda\left(1+\frac{x}{\lambda'} \right)}{2}e^{x\left(\frac{\lambda}{\lambda'} -1\right)} \\
&\quad + \frac{\lambda\left(1-\frac{x}{\lambda'} \right)}{2}e^{-x\left(\frac{\lambda}{\lambda'} -1\right)} - \lambda \\
&= 2 \lambda \sinh^2\left( \frac{x}{2} \left(\frac{\lambda}{\lambda'}-1\right) \right) \\
&\quad + x\frac{\lambda}{\lambda'} \sinh \left( x \left(\frac{\lambda}{\lambda'}-1\right) \right).
\end{align*}
This proves (ii). 
For (iii), by Stirling's approximation, the maximal binomial coefficient $\binom{m}{k}$ is $O(2^m/\sqrt{m})$ which is $o(e^{m})$. Therefore,  
\begin{align*}
    2\left(1+g_{m,\lambda}(t \lambda)\right) &= e^{-\lambda t}(1+t)^m + e^{\lambda t}(1-t)^m\\
    &= \left(1 - \lambda t + \lambda^2 t^2/2 + o(t^2 \lambda^2)\right) \\
    &\quad \times \left(1 + mt + m(m-1)t^2/2 + o(t^2 e^m) \right) \\
    &\quad + \left(1 + \lambda t + \lambda^2 t^2/2 + o(t^2 \lambda^2)\right) \\
    &\quad \times \left(1 - mt + m(m-1)t^2/2 + o(t^2e^m) \right) \\
    &= 2 + t^2 \left( m(m-1) - 2 \lambda m + \lambda^2 \right) \\
    &\quad +o(t^2 (e^m + \lambda^2) ). 
\end{align*}

Notice that for $t \in  [-c\lambda, c \lambda]$ we have 
\begin{align*}
    2(1+g_{m,\lambda}(|t|)) & = e^{-|t|}(1+|t|/\lambda)^m + e^{|t|}(1-|t|/\lambda)^m \\
    & \geq (1-c)^m. 
\end{align*}
Therefore, $g_{m,\lambda}(|t|) \geq C :=(1-c)^m/2 - 1$, where $C > -1$ because $c<1$. Because $g_{m,\lambda}(t)$ is symmetric around $t=0$, we have $g_{m,\lambda}(t) \geq C$. This proves (iv). \qed

The following lemma provides identities of sums of the functions $h_{m,\lambda}$ and $g_{m,\lambda}$. 
\begin{lem}
\label{lem:gm_sums}
\begin{enumerate}
The following holds for $t,s \in \reals$.
\item[(i)] 
\begin{align}
H(t,s) := H(t,s;\lambda) & := \sum_{m=0}^\infty \Po_{\lambda}(m)h_{m,\lambda}(t)h_{m,\lambda}(s) \nonumber \\
& =  e^{ts/\lambda}-1.
\label{eq:hm_identity}
\end{align}
\item[(ii)] 
\begin{align}
\label{eq:gm_second_moemnt}
G(t,s) &:= G(t,s;\lambda) :=  \sum_{m=0}^\infty \Po_{\lambda}(m)g_{m,\lambda}(t)g_{m,\lambda}(s) \nonumber \\
& = 2 \sinh^2 \left(\frac{ts}{2\lambda}\right).
\end{align}
\item[(iii)] As $u = t^2/\lambda \to 0$ and $\lambda \to 0$, 
\begin{align}
F(t) & := F(t;\lambda) := \sum_{m=0}^\infty \Po_{\lambda}(m)g_{m,\lambda}^4(t) \nonumber \\
& = u^4 \left( \frac{15}{4} + \frac{9}{\lambda} + \frac{1}{2 \lambda^2} \right) + o(u^4/\lambda).
\label{eq:gm_fourth_moemnt}
\end{align}
\end{enumerate}
\end{lem}
\subsection*{Proof of Lemma~\ref{lem:gm_sums}}
For (i), let $p=1+t/\lambda$ and $q=1+s/\lambda$. 
Then
\begin{align*}
H(t,s) &= \sum_{m=0}^{\infty} \frac{e^{-\lambda-t-s} (\lambda pq)^m}{m!} \\
&\quad - \sum_{m=0}^{\infty} \frac{e^{-\lambda-t} (\lambda p)^m}{m!} \\
&\quad - \sum_{m=0}^{\infty} \frac{e^{-\lambda-s} (\lambda q)^m}{m!} \\
&\quad + \sum_{m=0}^{\infty} \frac{e^{-\lambda} \lambda^m}{m!}. 
\end{align*}
By a standard expansion of the exponential function
\begin{align*}
H(t,s) &= e^{-\lambda-t-s} e^{\lambda pq} - e^{-\lambda-t} e^{\lambda p} \\
&\quad - e^{-\lambda-s} e^{\lambda q} + e^{-\lambda} e^{\lambda}. 
\end{align*}
Substituting $p$ and $q$ and simplifying leads to \eqref{eq:hm_identity}. (ii) follows similarly. For (iii), we first note that similar evaluations as in (i) give, for all $k=1,2,\ldots$ and $|t|\leq \lambda$,
\begin{align*}
    E_k(t) & := E_k(t;\lambda) \\
    &:= \sum_{m=0}^\infty (1+g_{m,\lambda}(t))^k \Po_\lambda(m)  \\
    &= e^{-\lambda}2^{-k} \sum_{j=0}^k \exp\{ \lambda(t(k-2j) \\
    &\quad + (1+\frac{t}{\lambda})^j (1-\frac{t}{\lambda})^{k-j})\}.
\end{align*}
Using the binomial expansion identity, we can express $F(t)$ using $\{E_k(t)\}_{k=1}^4$ as 
\begin{align*}
F(t) & = \sum_{k=0}^4 \binom{4}{k}(-1)^{4-k} E_k(t ; \lambda) \\
& = E_4(t) - 4E_3(t) + 6 E_2(t) - 4E_1(t) + E_0. 
\end{align*}
Set $u := t^2/\lambda$. For $r=1$, we have the following closed-form expressions: $E_0 = E_1 = 1$, $E_2(t) = \cosh(u)$, 
\begin{align*}
    4E_3(t) & = \cosh(u^2/t)(e^{3u} + 3e^{-u}), \\
    8 E_4(t) & = 4e^{-u^2/\lambda}\cosh(2u^2/t) + 3e^{u^2/\lambda-2u}+  e^{u^2/\lambda + 6u} \cosh(4u^2/t). 
\end{align*}
Putting it all together, we obtain:
\begin{align}
F(t) &= -3 + 6\cosh(u) -\cosh(u^2/\lambda)\,\bigl(e^{3u}+3e^{-u}\bigr) \label{eq:F_closed_form} \\
&\quad + \frac{3}{8}\,e^{\,u^2/t-2u} + \frac{1}{8}\,e^{u^2/t+6u}\cosh(4u^2/\lambda) \nonumber \\
&\quad + \frac{1}{2}e^{-u^2/t}\cosh(2u^2/\lambda).
\end{align}
Expanding in a Taylor series leads to
This leads to
\begin{align*}
F(t) &= u^4\Bigl(\frac{15}{4} + \frac{9}{\lambda} + \frac{1}{2\lambda^2}\Bigr) \\
&\quad + u^5\Bigl(6 + \frac{19}{\lambda}\Bigr) \\
&\quad + u^6\Bigl(\frac{57}{8} + \frac{41}{\lambda} + \frac{9}{\lambda^2}\Bigr) + O(u^7).
\end{align*}
This implies \eqref{eq:gm_fourth_moemnt}.
\qed

The following Lemma presents expansions used in the proof of Theorem~\ref{thm:lower_bound}. 
\begin{lem} For $\lambda > 0$ and $|t| < \lambda$, 
$-1 < g_{m,\lambda}(t)$ by Lemma~\ref{lem:gm_properties}-(iv). Thus, 
\begin{align*}
 Q_k(t ; \lambda) := \sum_{m=0}^\infty \left|\log(1+g_{m,\lambda}(t))\right|^k \Po_\lambda(m),
\end{align*} 
is defined for any $k=1,2,\ldots$. The following holds. Set $u := t^2/\lambda$. If $u \to 0$ while $\lambda = O(1)$, then 
\begin{enumerate}
\item[(i)]  $Q_1(t; \lambda) = o(u)$. 
\item[(ii)] $Q_2(t ; \lambda) = \frac{1}{2} u^2 (1 + o(1))$.
\item[(iii)] $Q_4(t ; \lambda) = u^4 \left(\frac{15}{4} + \frac{9}{\lambda} + \frac{1}{2 \lambda^2} + o(1)\right)$.
\end{enumerate}
\label{lem:gm_log_gm}
\end{lem}

\subsection*{Proof of Lemma~\ref{lem:gm_log_gm}}
By the inequality $\log(1+x) \leq x$ and Lemma~\ref{lem:gm_properties}-(i), we get
\[
Q_1(t; \lambda) \leq  \sum_{m=0}^\infty g_{m,\lambda}(t) \Po_\lambda(m) = 0.
\]
Additionally, by Lemma~\ref{lem:gm_properties}-(iii),
we have 
\begin{align*}
g_{m,\lambda}(t) &= \frac{u}{2} \left( m(m-1)/\lambda - 2 + \lambda \right) \\
&\quad + o (u(e^m + \lambda^3)). 
\end{align*}
Because the Poisson moment generation function exists, we can write
\begin{align*}
\sum_{m=0}^\infty o(u(e^m + \lambda^3)) \Po_{\lambda}(m) &= o( (u e^{\lambda(e-1)}+u\lambda^3)) \\
&= o(u). 
\end{align*}
From the Taylor expansion $\log(1+z) = z (1 + O(z))$, it follows that 
\begin{align*}
\bigl(\log(1+g_{m,\lambda}(t))\bigr)^k
&= g_{m,\lambda}(t)^k \left(1 + O(g_{m,\lambda}(t)) \right) \\
&= g_{m,\lambda}(t)^k + O(g_{m,\lambda}^{k+1}(t)).
\end{align*}
The identity $\sum_{m=0}^\infty g_{m,\lambda}(t) \Po_{\lambda}(m)=0$ implies (i). Putting the above in the expressions for $Q_2(t;\lambda)$, we get
\begin{align*}
Q_2(t;\lambda)
&= \sum_{m=0}^\infty g_{m,\lambda}(t)^2 \Po_\lambda(m)
   + O(g_{m,\lambda}(t)^3) \\
&= G(t,t; \lambda) + o(u^3),
\end{align*}
where $G(t,s;\lambda)$ is provided in Lemma~\ref{lem:gm_sums}-(ii), and
\begin{align*}
Q_4(t;\lambda)
&= \sum_{m=0}^\infty g_{m,\lambda}(t)^4 \Po_\lambda(m)
   + O(g_{m,\lambda}(t)^4) \\
&= F(t;\lambda) + o(u^4),
\end{align*}
where $F(t;\lambda)$ is provided in Lemma~\ref{lem:gm_sums}-(iii). From Lemma~\ref{lem:gm_sums} we also get the expansion of $F(t;\lambda)$ in powers of $u = t^2/\lambda$, leading to 
\[
Q_4(t;\lambda)
= u^4\Bigl(\frac{15}{4} + \frac{9}{\lambda} + \frac{1}{2\lambda^2}\Bigr)
  + O(u^5).
\]
A similar expansion leads to the expression of $Q_2(t;\lambda)$, so (ii) and (iii) follows.
\qed

\section*{Acknowledgment}
The author would like to thank David Donoho for discussing the problem formulation, Reut Levi for comments and suggestions concerning an earlier version of this manuscript, and several anonymous reviewers whose feedback significantly improved the scope of the results and the quality of this paper.

\bibliography{multinomials,IEEEfull}

\end{document}